\documentclass{amsart}

\usepackage{amssymb}
\usepackage[OT2,T1]{fontenc}
\usepackage{eepic}
\usepackage{color}

\allowdisplaybreaks

\numberwithin{equation}{section}

\DeclareSymbolFont{cyrletters}{OT2}{wncyr}{m}{n}
\DeclareMathSymbol{\Sha}{\mathalpha}{cyrletters}{"58}

\newtheorem{theorem}{Theorem}[section]
\newtheorem{lemma}[theorem]{Lemma}

\newtheorem{remark}[theorem]{Remark}

\newtheorem{TheoA}{Theorem A}
\newtheorem{TheoB}{Theorem B}
\newtheorem{TheoC}{Theorem C}

\newtheorem{Cuculescutheo}{Cuculescu's construction \cite{Cu}}

\newtheorem{CZDecomposition}{Calder\'on-Zygmund decomposition \cite{Pa1}}
\newtheorem{Atomic}{Atoms and John-Nirenberg inequality \cite{Pe2,HM}}

\newcommand{\Z}{\mathbb{Z}}
\newcommand{\R}{\mathbb{R}}
\newcommand{\C}{\mathbb{C}}

\newcommand{\summ}{\sum\nolimits}

\def\G{\mathrm{G}}
\def\1{\mathbf{1}}

\def\Q{\mathcal{Q}}
\def\M{\mathcal{M}}
\def\A{\mathcal{A}}

\def\LT{\mathsf{LT}}
\def\UT{\mathsf{UT}}

\newcommand{\dem}{\noindent {\bf Proof. }}
\newcommand{\demAia}{\noindent {\bf Proof of Theorem Ai) --- Perfect dyadic CZO's. }}
\newcommand{\demAib}{\noindent {\bf Proof of Theorem Ai) --- Haar shift operators. }}
\newcommand{\demAii}{\noindent {\bf Proof of Theorem Aii). }}
\newcommand{\demB}{\noindent {\bf Proof of Theorem B. }}
\newcommand{\demCi}{\noindent {\bf Proof of Theorem C --- Weak type inequalities. }}
\newcommand{\demCii}{\noindent {\bf Proof of Theorem C --- $\mathrm{H}_p/L_p$ type inequalities. }}

\newcommand{\fin}{\hspace*{\fill} $\square$ \vskip0.2cm}

\def\mean{- \hskip-10.6pt \int}

\addtocounter{tocdepth}{-1}

\begin{document}

\title[CZO's with matrix-valued kernels]{Calder\'on-Zygmund operators \\ associated to matrix-valued kernels}

\author[Hong, L\'opez-S\'anchez, Martell and Parcet]
{Guixiang Hong, Luis Daniel L\'opez-S\'anchez, \\ Jos\'e Mar\'ia Martell and Javier Parcet}

\null

\vskip-40pt

\null

\maketitle

\vskip-30pt

\null

\begin{abstract} 
Calder\'on-Zygmund operators with noncommuting kernels may fail to be $L_p$-bounded for $p \neq 2$, even for kernels with good size and smoothness properties. Matrix-valued paraproducts, Fourier multipliers on group vNa's or noncommutative martingale transforms are frameworks where we find such difficulties. We obtain weak type estimates for perfect dyadic CZO's and cancellative Haar shifts associated to noncommuting kernels in terms of a row/column decomposition of the function. Arbitrary CZO's satisfy $\mathrm{H}_1 \to L_1$ type estimates. In conjunction with $L_\infty \to \mathrm{BMO}$, we get certain row/column $L_p$ estimates. Our approach also applies to noncommutative paraproducts or martingale transforms with noncommuting symbols/coefficients. Our results complement recent results of Junge, Mei, Parcet and Randrianantoanina.    
\end{abstract}

\addtolength{\parskip}{+1ex}

\section*{Introduction}

A \emph{semicommutative} CZO has the formal expression $$Tf(x) \, \sim \, \int_{\R^n} k(x,y) (f(y)) \, dy,$$ where the kernel acts linearly on the matrix-valued function $f = (f_{ij})$ and satisfies standard size/smoothness Calder\'on-Zygmund type conditions. This is the operator model for quite a number of problems which have attracted some attention in recent years, including matrix-valued paraproducts, operator-valued Calder\'on-Zygmund theory or Fourier multipliers on group von Neumann algebras, see \cite{JM,JMP1,M1,NPTV,Pa1} and the references therein. To be more precise, let $\mathcal{B}(\ell_2)$ stand for the matrix algebra of bounded linear operators on $\ell_2$. Consider the algebra formed by essentially bounded functions $f: \R^n \to \mathcal{B}(\ell_2)$. Its weak operator closure is a von Neumann algebra $\mathcal{A}$ and as such we may construct  noncommutative $L_p$ spaces over it. Let us highlight a few significant examples:

\begin{itemize}
\item \textbf{Scalar kernels.} $k(x,y) \in \C$ and $$\hskip-32pt k(x,y)(f(y)) \ = \ \Big( k(x,y) f_{ij}(y) \Big).$$ 

\vskip3pt

\item \textbf{Schur product actions.} $k(x,y) \in \mathcal{B}(\ell_2)$ and $$\hskip-26pt k(x,y)(f(y)) \ = \ \Big( k_{ij}(x,y) f_{ij}(y) \Big).$$

\vskip3pt

\item \textbf{Fully noncommutative model.} $k(x,y) \in \mathcal{B}(\ell_2) \bar\otimes \mathcal{B}(\ell_2)$ and 
$$\hskip28pt k(x,y)(f(y)) \ = \ \Big( \summ_m \mathrm{tr} \big( k_m''(y) f(y) \big) k_m'(x)_{ij}\Big).$$

\vskip3pt

\item \textbf{Partial traces, noncommuting kernels.} $k(x,y) \in \mathcal{B}(\ell_2)$ and $$\hskip14pt k(x,y) (f(y)) \ = \ \left\{ \begin{array}{c} \Big( \sum_s k_{is}(x,y) f_{sj}(y) \Big), \\ [7pt] \Big( \sum_s f_{is}(y) k_{sj}(x,y) \Big). \end{array} \right.$$
\end{itemize}

Scalar kernels required in \cite{Pa1} a matrix-valued Calder\'on-Zygmund decomposition in terms of noncommutative martingales and a pseudo-localization principle to control the tails of $T \hskip-1pt f$ in the $L_2$-metric. Hilbert space valued kernels were later considered in \cite{MP}, see also \cite{MeiMAMS,PX1,Rad} for previous related results. The second case refers to the Schur matrix product $k(x,y) \bullet f(y)$, considered for the first time in \cite{JMP1} to analyze cross product extensions of classical CZO's. It is instrumental for H\"ormander-Mihlin type theorems on Fourier multipliers associated to discrete groups and for Schur multipliers with a Calder\'on-Zygmund behavior \cite{JMP1,JMP2}. In the fully noncommutative model, we approximate $k(x,y)$ by a sum of elementary tensors $\sum_m k_m'(x) \otimes k_m''(y)$ and the action is given by $$Tf(x) \ \sim \ \int_{\R^n} (id \otimes \mathrm{tr}) \, \Big[ k(x,y) \big( \mathbf{1} \otimes f(y) \big) \Big] \, dy.$$ In this case, we regard the space $L_p(\mathcal{A}) = L_p(\R^n; L_p(\mathcal{B}(\ell_2)))$ as a whole. In other words, the noncommutative nature of $L_p(\mathcal{A})$ predominates and the presence of a Euclidean subspace is ignored. That is what happens for purely noncommutative CZO's \cite{JMP3} and justifies the presence of $id \otimes \mathrm{tr}$, to integrate over the full algebra $\mathcal{A}$ and not just over the Euclidean part. The last case refers to matrix-valued kernels acting on $f$ by left/right multiplication, $k(x,y) f(y)$ and $f(y) k(x,y)$. Matrix-valued paraproducts are prominent examples \cite{K,M1,M2,NPTV,Po}. This is the only case in which the kernel does not commute with $f$, since the Schur product is abelian and we find $(id \otimes \mathrm{tr}) [ k(x,y) ( \mathbf{1} \otimes f(y) ) ] = (id \otimes \mathrm{tr}) [ ( \mathbf{1} \otimes f(y) ) k(x,y) ]$ by traciality. 

Our main goal is to obtain endpoint estimates for CZO's with noncommuting kernels, motivated by a recent estimate from \cite{JMP1} for semicommutative CZO's. If $k(x,y)$ acts linearly on $\mathcal{B}(\ell_2)$ and satisfies the H\"ormander smoothness condition in the norm of bounded linear maps on $\mathcal{B}(\ell_2)$, the content of \cite[Lemma 1.3]{JMP1} can be summarized as follows  
\begin{itemize}
\item If $T$ is $L_\infty(\mathcal{B}(\ell_2);L_2^r(\R^n))$-bounded, then $T: L_\infty(\mathcal{A}) \to \mathrm{BMO}_r(\mathcal{A})$,

\vskip1pt

\item If $T$ is $L_\infty(\mathcal{B}(\ell_2);L_2^c(\R^n))$-bounded, then $T: L_\infty(\mathcal{A}) \to \mathrm{BMO}_c(\mathcal{A})$.
\end{itemize}
Here, the $L_\infty(L_2^c)$-boundedness assumption refers to $$\Big\| \Big( \int_{\R^n} Tf(x)^* Tf(x) \, dx \Big)^\frac12 \Big\|_{\mathcal{B}(\ell_2)} \lesssim \ \Big\| \Big( \int_{\R^n} f(x)^* f(x) \, dx \Big)^\frac12 \Big\|_{\mathcal{B}(\ell_2)},$$ while the column-BMO norm of a matrix-valued function $g$ is given by $$\sup_{Q \ \mathrm{cube}} \Big\| \Big( \mean_Q \big( g(x) - g_Q \big)^* \big( g(x) - g_Q \big) \, dx \Big)^\frac12 \Big\|_{\mathcal{B}(\ell_2)}.$$ Taking adjoints ---so that the $*$ switches everywhere from left to right--- we find $L_\infty(L_2^r)$-boundedness and the row-BMO norm. The noncommutative BMO space $\mathrm{BMO}(\mathcal{A}) = \mathrm{BMO}_r(\mathcal{A}) \cap \mathrm{BMO}_c(\mathcal{A})$ was introduced in \cite{PX1}. According to \cite{Musat} it has the expected interpolation behavior in the $L_p$ scale. Thus, standard interpolation and duality arguments show that $T: L_p(\mathcal{A}) \to L_p(\mathcal{A})$ for $1 < p < \infty$ provided the kernel is smooth enough in both variables and $T$ is a normal self-adjoint map satisfying the $L_\infty(L_2^r)$ and $L_\infty(L_2^c)$ boundedness assumptions. In other words, the row/column boundedness conditions essentially play the role of the $L_2$-boundedness assumption in classical Calder\'on-Zygmund theory. 

Although this certainly works for non-scalar kernels ---Schur product actions were used e.g. in \cite[Theorem B]{JMP1}--- the boundedness assumptions impose nearly commuting conditions on the kernel which are too strong for CZO's associated to noncommuting kernels. Namely, given $k: \R^{2n} \setminus \Delta \to \mathcal{B}(\ell_2)$ smooth and given $x \notin \mathrm{supp}_{\R^n} f$, let us set formally the row/column CZO's $$T_cf(x) = \int_{\R^n} k(x,y) f(y) \, dy \quad \mbox{and} \quad T_rf(x) = \int_{\R^n} f(y) k(x,y) \, dy.$$ 
It is not difficult to construct noncommuting kernels with
\begin{itemize}
\item[i)] $T_r$ and $T_c$ are $L_2(\mathcal{A})$-bounded,

\vskip1pt

\item[ii)] $T_r$ and $T_c$ are not $L_p(\mathcal{A})$-bounded for $1 < p \neq 2 < \infty$, 
\end{itemize}
see e.g. \cite[Section 6.1]{Pa1} for specific examples. Therefore, the $L_\infty(L_2^r)$ and $L_\infty(L_2^c)$ boundedness assumption is in general too restrictive when kernel and function do not commute. Assume for what follows that $T_r$ and $T_c$ are $L_2(\mathcal{A})$-bounded. We are interested in weakened forms of $L_p$ boundedness and endpoint estimates for these CZO's. A \emph{dyadic noncommuting} CZO will be a $L_2(\mathcal{A})$-bounded pair $(T_r,T_c)$ associated to a noncommuting kernel satisfying one of the following conditions:
\begin{itemize}
\item[\textbf{a)}] \textbf{Perfect dyadic kernels} $$\big\| k(x,y) - k(z,y) \big\|_{\mathcal{B}(\ell_2)} + \big\| k(y,x) - k(y,z) \big\|_{\mathcal{B}(\ell_2)} = 0$$ whenever $x,z \in Q$ and $y \in R$  for some disjoint dyadic cubes $Q,R$.

\vskip5pt

\item[\textbf{b)}] \textbf{Cancellative Haar shift operators} $$k(x,y) = \sum_{Q \ \mathrm{dyadic}} \ \sum_{\begin{subarray}{c} R,S \ \mathrm{dyadic} \, \subset \, Q \\ \ell(R) = 2^{-r} \ell(Q) \\ \ell(S) = 2^{-s} \ell(S) \end{subarray}} \alpha_{RS}^Q h_R(x) h_S(y),$$ for some fixed $r,s \in \mathbb{Z}_+$ where the $\alpha_{RS}^Q \in \mathcal{B}(\ell_2)$ with $\|\alpha_{RS}^Q\|_{\mathcal{B}(\ell_2)} \le \frac{\sqrt{|R||S|}}{|Q|}$. Here $h_Q$ refers to any of the $2^n-1$ Haar functions related to the cube $Q$. 
\end{itemize}
Perfect dyadic kernels were introduced in \cite{AHMTT} and include Haar multipliers, as well as paraproducts and their adjoints. If $J_-$ and $J_+$ denote the left/right halves of a dyadic interval in $\R$, the standard model for Haar shifts is the dyadic Hilbert transform with kernel $\sum_J (h_{J_-}(y) - h_{J_+}(y)) h_J(x)$. It appeared after Petermichl's crucial result \cite{Pe}, showing the classical Hilbert transform as a certain average of dyadic Hilbert transforms. Hyt\"onen's representation theorem \cite{Hy} extends this result to arbitrary CZO's. We will write \emph{generic noncommuting} CZO for $L_2(\mathcal{A})$-bounded pairs $(T_r,T_c)$  with a noncommuting kernel satisfying the standard smoothness. Our first significant result is the following.     

\begin{TheoA} The following inequalities hold$\, :$
\begin{itemize}
\item[i)] \emph{Dyadic noncommuting CZO's}. Given $f \in L_1(\mathcal{A})$ $$\inf_{f = f_r + f_c} \big\|T_r f_r\big\|_{1,\infty} + \big\|T_c f_c \big\|_{1,\infty} \, \lesssim \, \|f\|_1.$$

\vskip3pt

\item[ii)] \emph{Generic noncommuting CZO's}. Given $f \in \mathrm{H}_1(\mathcal{A})$ $$\inf_{f = f_r + f_c} \big\|T_r f_r\big\|_1 + \big\|T_c f_c \big\|_1 \, \lesssim \, \|f\|_{\mathrm{H}_1(\mathcal{A})}.$$
\end{itemize}
\end{TheoA}

The noncommutative forms of $L_{1,\infty}$ and the Hardy space $\mathrm{H}_1$ are well-known in the subject. Nevertheless, they will also be properly defined in the body of the paper. Our main result is the inequality given in Theorem A i) and their noncommutative generalizations in Theorem C below. As we shall explain in the Appendix, the left/right modular nature of $T_r/T_c$ is essential for the weak type $(1,1)$ estimates, see also Remark \ref{Rem-Modularity}. The following result easily follows from Theorem A by interpolation/duality and it can also be derived from \cite{JMP1}. Nevertheless, it is worth mentioning the $L_p$ inequalities that we find. 

\begin{TheoB}
The following inequalities hold for generic noncommuting CZO's$\, :$
\begin{itemize}
\item[i)] If $1 < p < 2$ and $f \in L_p(\mathcal{A})$ $$\inf_{f = f_r + f_c} \big\|T_r f_r\big\|_p + \big\|T_c f_c \big\|_p \, \lesssim \, \|f\|_p.$$

\vskip3pt

\item[ii)] If $2 < p < \infty$ and $f \in L_p(\mathcal{A})$ $$\big\|T_r f\big\|_{\mathrm{H}_p^r(\mathcal{A})} + \big\|T_c f \big\|_{\mathrm{H}_p^c(\mathcal{A})} \, \lesssim \, \|f\|_p.$$

\vskip3pt

\item[iii)] Given $f \in L_\infty(\mathcal{A})$, we also have $\|T_rf\|_{\mathrm{BMO}_r(\mathcal{A})} + \|T_cf\|_{\mathrm{BMO}_c(\mathcal{A})} \lesssim \|f\|_\infty$.
\end{itemize}
\end{TheoB}

Theorems A and B also hold for other operator-valued functions, replacing $\mathcal{B}(\ell_2)$ by any semifinite von Neumann algebra $\M$. Our proof will be written in this framework. Let us now consider a weak-$*$ dense filtration $\Sigma_\A = (\A_n)_{n \ge 1}$ of von Neumann subalgebras of an arbitrary semifinite von Neumann algebra $\A$. In the following result, we will consider two kind of operators in $L_p(\A)$: 
\begin{itemize} 
\item[\textbf{a)}] \textbf{Noncommuting martingale transforms} $$M_\xi^r f = \sum_{k \ge 1} \Delta_k(f) \xi_{k-1} \quad \mbox{and} \quad M_\xi^c f = \sum_{k \ge 1} \xi_{k-1} \Delta_k(f).$$ 

\vskip3pt

\item[\textbf{b)}] \textbf{Paraproducts with noncommuting symbol} $$\Pi_\rho^r(f) = \sum_{k \ge 1} \mathsf{E}_{k-1}(f) \Delta_k(\rho)  \quad \mbox{and} \quad \Pi_\rho^c(f) = \sum_{k \ge 1} \Delta_k(\rho) \mathsf{E}_{k-1}(f).$$
\end{itemize}
Here $\Delta_k$ denotes the martingale difference operator $\mathsf{E}_k - \mathsf{E}_{k-1}$ and $\xi_{k} \in \A_{k}$ is an adapted sequence. Of course, the symbols $\xi$ and $\rho$ do not necessarily commute with the function. Randrianantoanina considered in \cite{Rad} noncommutative martingale transforms with commuting coefficients. As for paraproducts with noncommuting symbols, Mei studied the $L_p$-boundedness for $p>2$ and regular filtrations in \cite{M1} and also analyzed in \cite{M2} the case $p<2$ in the dyadic matrix-valued case under a strong BMO condition of the symbol. Our theorem below goes beyond these results, see also \cite{MP} for related results.

\begin{TheoC}
Consider the pairs$\,:$
\begin{itemize}
\item[i)] \emph{Martingale transforms} $(M_\xi^r, M_\xi^c)$, with $\sup_k \|\xi_k\|_\M < \infty$.

\item[ii)] \emph{Martingale paraproducts} $(\Pi_\rho^r, \Pi_\rho^c)$, with $\Pi_\rho^{r/c}$ $L_2(\A)$-bounded.
\end{itemize}
If $\Sigma_\A$ is regular, we obtain weak type $(1,1)$ inequalities like in Theorem \emph{Ai)} for martingale transforms and paraproducts . The estimates in Theorems \emph{Aii)} and \emph{B} also hold for both families and for arbitrary filtrations $\Sigma_\A$. Moreover, the martingale paraproducts $\Pi_\rho^r$ and $\Pi_\rho^c$ are $L_p$-bounded for $2 < p < \infty$ and $L_\infty \to \mathrm{BMO}$. 
\end{TheoC}

In the case of martingale transforms, there are also examples of noncommuting kernels failing $L_p$-boundedness for $p \neq 2$. Hence, our results recover those in \cite{Rad,Rad2} and are in some sense sharp, providing appropriate substitutes for noncommuting coefficients. Our result for paraproducts goes beyond \cite[Theorem 1.2]{M1} in two aspects. First, our estimates for $p > 2$ hold for arbitrary martingales, not just for regular ones. Second, we give a partial answer to Mei's question in \cite{M1} after the proof of Theorem 1.2 for the case $p < 2$ and also for the weak type $(1,1)$ estimates. The paper is organized following the order in the Introduction. We include an Appendix at the end with further comments and open problems. Along the paper we shall assume some familiarity with basic notions from noncommutative integration. The content of \cite[Section 1]{Pa1} is enough for our purposes, more can be found in \cite{KR,PX2,Ta}.

\section{Calder\'on-Zygmund decomposition}

Let $\M$ be a semifinite von Neumann algebra equipped with a normal semifinite faithful trace $\tau$. Consider the algebra of essentially bounded functions $\R^n \to \M$ equipped with the n.s.f. trace $$\varphi(f) = \int_{\R^n} \tau(f(x)) \, dx.$$ Its weak-operator closure is a von Neumann algebra $\A$. If $1 \le p \le \infty$, we write $L_p(\M)$ and $L_p(\A)$ for the noncommutative $L_p$ spaces associated to the pairs $(\M,\tau)$ and $(\A,\varphi)$. The lattices of projections are written $\M_\pi$ and $\A_{\pi}$, while $\1_\M$ and $\1_{\A}$ stand for the unit elements. The set of dyadic cubes in $\R^n$ is denoted by $\Q$ and we use $\Q_k$ for the $k$-th generation, formed by cubes $Q$ with side length $\ell(Q) = 2^{-k}$. If $f: \R^n \to \M$ is integrable on $Q \in \Q$, we set the average $$f_Q = \frac{1}{|Q|} \int_Q f(y) \, dy.$$ Let us write $(\mathsf{E}_k)_{k \in \Z}$ for the family of conditional expectations associated to the classical dyadic filtration on $\R^n$. $\mathsf{E}_k$ will also stand for the tensor product $\mathsf{E}_k \otimes id_\M$ acting on $\A$. If $1 \le p \le \infty$ and $f \in L_p(\A)$ $$\begin{array}{rclcl}
\mathsf{E}_k(f) & = & f_k & = & \displaystyle \sum_{Q \in \Q_k}^{\null} f_Q 1_Q, \\ [15pt] \Delta_k(f) & = & df_k & = & \displaystyle \sum_{Q \in \Q_k} \big( f_Q - f_{\widehat{Q}} \big) 1_Q, \end{array}$$ where $\widehat{Q}$ denotes the dyadic parent of $Q$. We will write $(\A_k)_{k \in \Z}$ for the filtration $\A_k = \mathsf{E}_k(\A)$. The noncommutative weak $L_1$-space, denoted by $L_{1,\infty}(\mathcal{A})$, is the set of all $\varphi$-measurable operators $f$ for which $\left\|f\right\|_{1,\infty} = \sup_{\lambda > 0} \, \lambda \hskip1pt \varphi \{ |f| > \lambda \} < \infty$, see \cite{FK} for a more in depth discussion. In this case, we write $\varphi \{ |f| > \lambda \}$ to denote the trace of the spectral projection of $|f|$ associated to the interval $(\lambda,\infty)$. We find this terminology more intuitive, since it is reminiscent of the classical one. The space $L_{1,\infty}(\A)$ is a quasi-Banach space and satisfies the quasi-triangle inequality below which will be used with no further reference $$\lambda \, \varphi \Big\{ |f_1+f_2| > \lambda \Big\} \le \lambda \, \varphi \Big\{ |f_1| > \lambda/2 \Big\} + \lambda \, \varphi \Big\{ |f_2| > \lambda/2 \Big\}.$$ Let us consider the dense subspace $$\A_{c,+} = L_1(\A) \cap \Big\{ f: \R^n \to \M \, \big| \ f \in \A_+, \ \mathrm{supp}_{\R^n} \hskip1pt f \ \ \mathrm{is \ compact} \Big\} \subset L_1^+(\A).$$ Here 
$\mathrm{supp}_{\R^n}$ means the support of $f$ as a vector-valued function in $\R^n$. In other words, we have $\mathrm{supp}_{\R^n} \hskip1pt f = \mathrm{supp} \hskip1pt \|f\|_\M$. We employ this terminology to distinguish from $\mathrm{supp} \, f$, the support of $f$ as an operator in $\A$. Any function $f \in \A_{c,+}$ gives rise to a martingale $(f_k)_{k \in \Z}$ with respect to the dyadic filtration. Moreover, it is clear that given $f \in \A_{c,+}$ and $\lambda > 0$, there must exist $m_\lambda(f) \in \Z$ so that $0 \le f_k \le \lambda$ for all $k \le m_\lambda(f)$. The noncommutative analogue of the weak type $(1,1)$ boundedness of Doob's maximal function is due to Cuculescu. Here we state it in the context of operator-valued functions from $\A$. 

\begin{Cuculescutheo}
Let $f \in \A_{c,+}$ and consider the corresponding martingale $(f_k)_{k \in \Z}$ relative to the filtration $(\mathcal{A}_k)_{k \in \Z}$. Given $\lambda \in \R_+$, there exists a decreasing sequence of projections $(q_k(\lambda))_{k \in \Z}$ in $\mathcal{A}$ satisfying 
\begin{itemize}
\item[i)] $q_k(\lambda)$ commutes with $q_{k-1}(\lambda) f_k q_{k-1}(\lambda)$ for each $k$,
\item[ii)] $q_k(\lambda)$ belongs to $\mathcal{A}_k$ for each $k$ and $q_k(\lambda) f_k q_k(\lambda) \le \lambda \hskip1pt q_k(\lambda)$,
\item[iii)] The following estimate holds $$\varphi \Big( \mathbf{1}_\A - \bigwedge_{k \in \Z} q_k(\lambda) \Big) \ \le \
\frac{1}{\lambda} \hskip1pt \sup_{k \in \Z} \|f_k\|_1 \ = \ \frac{1}{\lambda} \hskip1pt \|f\|_1.$$
\end{itemize}
Explicitly, take $q_k(\lambda) = \chi_{(0,\lambda]}(q_{k-1}(\lambda) f_k q_{k-1}(\lambda))$ with $q_k(\lambda) = \1_\A$ for $k \le m_\lambda(f)$.
\end{Cuculescutheo}

Given $f \in \A_{c,+}$, consider the Cuculescu's sequence $(q_k(\lambda))_{k \in \Z}$ associated to $(f,\lambda)$ for a given $\lambda >0$. Since $\lambda$ will be fixed most of the time, we will shorten the notation by $q_k$ and only write $q_k(\lambda)$ when needed. Define the sequence $(p_k)_{k \in \Z}$ of disjoint projections $p_k =
q_{k-1}-q_k$, so that $$\sum_{k \in \Z} p_k = \mathbf{1}_\A - q \quad \mbox{with} \quad q = \bigwedge_{k \in
\Z} q_k.$$

\begin{CZDecomposition}
Given $f \in \A_{c,+}$ and $\lambda > 0$, we may decompose $f = g_d + g_\mathit{off} + b_d + b_\mathit{off}$ as the sum of four operators defined in terms of the Cuculescu's construction as follows
\begin{eqnarray*}
g_d & = & qfq + \sum_{k \in \Z} p_k f_k p_k, \\
b_d & = & \sum_{k \in \Z} p_k \hskip1pt (f - f_k) \hskip1pt p_k, \\
b_{\mathit{off}} & = & \sum_{i \neq j} p_i (f-f_{i \vee j}) p_j, \\
g_\mathit{off} & = & \sum_{i \neq j} p_i f_{i \vee j} p_j \ + \ q f (\mathbf{1}_\A - q) + (\mathbf{1}_\A - q) f q.
\end{eqnarray*}
Moreover, we have the diagonal estimates $$\Big\| qfq + \sum_{k \in \Z} p_k f_k p_k \Big\|_2^2 \le 2^n \lambda \, \|f\|_1 \quad \mbox{and} \quad \sum_{k \in \Z} \big\| p_k (f-f_k) p_k \big\|_1 \le 2 \, \|f\|_1.$$ The expression below for $g_{\mathit{off}}$ will be also instrumental $$g_{\mathit{off}} = \sum_{s=1}^\infty \sum_{k=m_\lambda+1}^\infty p_k df_{k+s} q_{k+s-1} + q_{k+s-1} df_{k+s} p_k = \sum_{s=1}^\infty \sum_{k=m_\lambda+1}^\infty g_{k,s} = \sum_{s=1}^\infty g_{(s)}.$$
\end{CZDecomposition}

\section{Proof of Theorems A and B}

The key result of this paper is Theorem A, since the remaining theorems follow from it or by using analog ideas. We begin with the proof of the weak type estimates for perfect dyadic CZO's and then make the necessary adjustments to make it work for Haar shift operators. The proof of Theorem Aii) will require to recall some recent results on square function and atomic Hardy spaces. 

\subsection{Perfect dyadic CZO's} 

To the best of our knowledge, the notion of perfect dyadic Calder\'on-Zygmund operator was rigorously defined for the first time in \cite{AHMTT} by Auscher, Hofmann, Muscalu, Tao and Thiele. Accordingly, we define a \emph{perfect dyadic \emph{CZO} with noncommuting kernel} as a pair $(T_r,T_c)$ formally given by 
\begin{eqnarray*}
T_rf(x) & \sim & \int_{\R^n} k(x,y) f(y) \, dy, \\
T_cf(x) & \sim & \int_{\R^n} f(y) k(x,y) \, dy,
\end{eqnarray*}
with an $\M$-valued kernel satisfying the perfect dyadic conditions $$\big\| k(x,y) - k(z,y) \big\|_\M + \big\| k(y,x) - k(y,z) \big\|_{\M} = 0$$ whenever $x,z \in Q$ and $y \in R$  for some disjoint dyadic cubes $Q,R$. Alternatively, we may think of perfect dyadic kernels $k: \R^{2n} \setminus \Delta \to \M$ as those which are constant on $2n$-cubes of the form $Q \times R$, where $Q,R$ are distinct dyadic cubes in $\R^n$ with the same side length and sharing the same dyadic parent. Classical perfect dyadic CZO's include Haar multipliers/martingale transforms and dyadic paraproducts. In other words, operators of the following form 
\begin{eqnarray*}
H_\xi f(x) & = & \int_{\R^n} \Big( \sum_{Q \in \Q} \frac{\xi(\widehat{Q})}{|Q|} \, 1_Q(x)(1_Q - 2^{-n}1_{\widehat{Q}})(y) \Big) f(y) \, dy, \\ \Pi_\rho f(x) & = & \int_{\R^n} \Big( \sum_{Q \in \Q} \frac{1}{|Q|} (\rho_Q - \rho_{\widehat{Q}}) 1_Q(x) 2^{-n}1_{\widehat{Q}}(y) \Big) f(y) \, dy,
\end{eqnarray*}
with $\sup_Q |\xi(Q)| < \infty$ and $\rho: \R^n \to \C$ in dyadic BMO. Adjoints of paraproducts are also perfect dyadic. In the noncommuting setting, the coefficients $\xi(Q)$ and the symbol $\rho$ become operators in $\M$ and an $\M$-valued function respectively which do not commute a priori with $f \in L_p(\A)$. Nevertheless, the perfect dyadic condition for the kernel is still satisfied in these cases.

\demAia Splitting $f$ as a sum of four positive operators and by density of $\A_{c,+}$ in the positive cone of $L_1(\A)$, we may clearly assume that $f \in \A_{c,+}$. A well-known lack of Cuculescu's construction is that we do not necessarily have $q_k(\lambda_1) \le q_k(\lambda_2)$ for $\lambda_1 \le \lambda_2$. This is typically solved restricting our attention to lacunary values for $\lambda$. Define $$\pi_{j,k} = \bigwedge_{s \ge j} q_k(2^s) - \bigwedge_{s \ge j-1} q_k(2^s) \quad \mbox{for} \quad j,k \in \Z.$$ We have $\sum_j \pi_{j,k} \stackrel{\mathrm{SOT}}{=} \1_\A - \psi_k$, where $$\psi_k = \bigwedge_{s \in \Z} q_k(2^s).$$ Observe that $\psi_k df_k = df_k \psi_k = 0$ for $k \in \Z$. Indeed, we have 
\begin{eqnarray*}
\|\psi_k df_k\|_\A & \le & \| \psi_k f_k^\frac12 \|_\A^{\null} \|f_k\|_\A^{\frac12} + \| \psi_{k} f_{k-1}^\frac12\|_\A^{\null {}} \|f_{k-1}\|_\A^{\frac12} \\ & = & \| \psi_k f_k \psi_k \|_\A^\frac12 \|f_k\|_\A^{\frac12} + \| \psi_{k} f_{k-1} \psi_k \|_\A^\frac12 \|f_{k-1}\|_\A^{\frac12} \ \le \ \lim_{s \to - \infty} 2^{1+\frac{s}{2}} \|f\|_\A^{\frac12}.
\end{eqnarray*}
In particular, we find $f = \sum_k (\1_\A - \psi_{k-1}) df_k (\1_\A - \psi_{k-1})$ and set $f = f_r + f_c$ with
\begin{eqnarray*}
f_r & = & \sum_{k \in \Z} \LT_{k-1}(df_k) \ = \ \sum_{k \in \Z} \Big( \sum_{i > j} \pi_{i,k-1} df_k \pi_{j,k-1} \Big), \\ f_c & = & \sum_{k \in \Z} \UT_{k-1}(df_k) \ = \ \sum_{k \in \Z} \Big( \sum_{i \le j} \pi_{i,k-1} df_k \pi_{j,k-1} \Big).
\end{eqnarray*}
This is the decomposition we will use for any perfect dyadic CZO. Given such an operator $T = (T_r,T_c)$ and $\lambda > 0$, the goal is to show that there exists an absolute constant $c_0$ so that $\lambda \varphi \{|T_rf_r| > \lambda\} + \lambda \varphi \{|T_cf_c| > \lambda\} \le c_0 \|f\|_1$ for any $f \in \A_{c,+}$ and any $\lambda > 0$. By symmetry in the argument, we will just prove the inequality for $T_cf_c$. Moreover, replacing $c_0$ by $2c_0$ we may also assume that $\lambda = 2^\ell$ for some $\ell \in \Z$. Having fixed the value of $\lambda$, we may consider the Calder\'on-Zygmund decomposition $f = g_d + g_{\mathit{off}} + b_d + b_{\mathit{off}}$ and set
$$\begin{array}{rclcrcl} g_d^c & = & \displaystyle \sum_{k \in \Z} \UT_{k-1} \big( \Delta_k (g_d) \big), & & g_{\mathit{off}}^c & = & \displaystyle \sum_{k \in \Z} \UT_{k-1} \big( \Delta_k (g_{\mathit{off}}) \big), \\ [15pt] b_d^c & = & \displaystyle \sum_{k \in \Z} \UT_{k-1} \big( \Delta_k (b_d) \big), & & b_{\mathit{off}}^c & = & \displaystyle \sum_{k \in \Z} \UT_{k-1} \big( \Delta_k (b_{\mathit{off}}) \big).
\end{array}$$  
By the quasi-triangle inequality it suffices to show $$\lambda \Big[ \varphi \Big\{ |T_c g_d^c| > \lambda \Big\} + \varphi \Big\{ |T_c b_d^c| > \lambda \Big\} + \varphi \Big\{ |T_c g_{\mathit{off}}^c| > \lambda \Big\} + \varphi \Big\{ |T_c b_{\mathit{off}}^c| > \lambda \Big\} \Big] \lesssim \|f\|_1.$$ The first term is first estimated by Chebychev's inequality in $\mathcal{A}$ $$\lambda \varphi \Big\{ |T_c g_d^c| > \lambda \Big\} \, \le \, \frac{1}{\lambda} \big\| T_c g_d^c \big\|_2^2 \, \lesssim \, \frac{1}{\lambda} \| g_d^c \|_2^2.$$ We use that $\UT_{k-1} \big( \Delta_k (g_d) \big)$ are in fact martingale differences, so that 
\begin{eqnarray*}
\frac{1}{\lambda} \| g_d^c \|_2^2 & = & \frac{1}{\lambda} \sum_{k \in \Z} \big\| \UT_{k-1} \big( \Delta_k (g_d)\big) \big\|_2^2 \ \le \ \frac{1}{\lambda} \sum_{k \in \Z} \| \Delta_k (g_d) \|_2^2 \\ & = & \frac{1}{\lambda} \Big\| \sum_{k \in \Z} \Delta_k (g_d) \Big\|_2^2 \ = \ \frac{1}{\lambda} \Big\| qfq + \sum_{k \in \Z} p_k f_k p_k \Big\|_2^2 \ \le \ 2^n \|f\|_1. 
\end{eqnarray*}
Indeed, the first inequality above follows from the fact that triangular truncations are contractive in $L_2(\A)$ while the last inequality arise from the diagonal estimates in the noncommutative CZ decomposition stated above. To handle the remaining terms, we introduce the projection $$\widehat{q} \, = \, \bigwedge_{s \ge \ell} q(2^s) \, = \, \bigwedge_{s \ge \ell} \bigwedge_{k \in \Z} q_k(2^s).$$ According to Cuculescu's construction, we find $$\varphi \big(\1_\A - \widehat{q} \hskip1pt \big) \, \le \, \sum_{s \ge \ell} \varphi \big( \1_\A - q(2^s) \big) \, \le \, \sum_{s \ge \ell} \frac{1}{2^s} \|f\|_1 \, = \, \frac{2}{\lambda} \|f\|_1.$$ This reduces our problem to show that $$\lambda \Big[ \varphi \Big\{ \big| T_c (b_d^c) \widehat{q} \hskip1pt \big| > \lambda \Big\} + \varphi \Big\{ \big| T_c (g_{\mathit{off}}^c) \widehat{q} \hskip1pt \big| > \lambda \Big\} + \varphi \Big\{ \big| T_c (b_{\mathit{off}}^c) \widehat{q} \hskip1pt \big| > \lambda \Big\} \Big] \ \lesssim \ \|f\|_1.$$ The perfect dyadic nature of $T_c$ comes now into scene. Indeed, we claim that the three terms $T_c (b_d^c) \widehat{q}, T_c (g_{\mathit{off}}^c) \widehat{q}, T_c (b_{\mathit{off}}^c) \widehat{q}$ vanish whenever $T_c$ is perfect dyadic. This will be enough to conclude the proof. If $Q_k(x)$ is the only cube in $\Q_k$ containing $x$, we find a.e. $x$
\begin{eqnarray*}
T_c(b_d^c)(x) \widehat{q}(x) & = & \sum_{k \in \Z} T_c \big( \UT_{k-1}(\Delta_k(b_d)) \big) (x) \, \widehat{q}(x) \\ & = & \sum_{k \in \Z} T_c \Big( \UT_{k-1}(\Delta_k(b_d)) 1_{Q_{k-1}(x)} \Big) (x) \, \widehat{q}(x) \\ & + & \sum_{k \in \Z} \sum_{\begin{subarray}{c} Q \in \Q_{k-1} \\ x \notin Q \end{subarray}} \Big( \int_Q k(x,y) \UT_{k-1}(\Delta_k(b_d)) (y) \, dy \Big) \, \widehat{q}(x). 
\end{eqnarray*}
The last term on the right vanishes since the term $\UT_{k-1}(\Delta_k(b_d))$ has mean $0$ in any $Q \in \Q_{k-1}$, so that we may replace $k(x,y)$ by $k(x,y) - k(x,c_Q)$, which is $0$ when $x \notin Q$ by the perfect dyadic cancellation of the kernel. On the other hand, if we define the projection $$\widehat{q}_{k-1} = \bigwedge_{s \ge \ell} q_{k-1}(2^s),$$ we see that $\widehat{q}(x) = \widehat{q}_{k-1}(x) \widehat{q}(x) = \widehat{q}_{k-1}(y) \widehat{q}(x)$ for any $y \in Q_{k-1}(x)$. This gives $$T_c(b_d^c)(x) \widehat{q}(x) \, = \, \summ_k T_c \Big( \UT_{k-1}(\Delta_k(b_d)) \widehat{q}_{k-1} 1_{Q_{k-1}(x)} \Big) (x) \, \widehat{q}(x).$$ The exact same argument applies for $g_{\mathit{off}}^c$ and $b_{\mathit{off}}^c$, so that it suffices to prove
\begin{eqnarray*}
\UT_{k-1}(\Delta_k(b_d)) \, \widehat{q}_{k-1} & = & 0, \\ \UT_{k-1}(\Delta_k(g_\mathit{off})) \, \widehat{q}_{k-1} & = & 0, \\ \UT_{k-1}(\Delta_k(b_\mathit{off})) \, \widehat{q}_{k-1} & = & 0,
\end{eqnarray*}
for all $k \in \Z$. In all these cases we will be using the following two key identities
\begin{itemize}
\item $\widehat{q}_{k-1} \pi_{i,k-1} = \pi_{j,k-1} \widehat{q}_{k-1} = 0$ for $i,j > \ell$ and $k \in \Z$,

\vskip3pt

\item $\pi_{i,k-1} p_{k-s} = p_{k-s} \pi_{j,k-1} = 0$ for $s \ge 1$, $i,j \le \ell$ and $k \in \Z$.
\end{itemize}
The proof is straightforward and left to the reader. It only requires to apply the monotonicity properties of $\bigwedge_{s \ge j} q_k(2^s)$, which increases in $j$ and decreases in $k$. If we apply the first identity to $\UT_{k-1}(\Delta_k(\gamma)) \, \widehat{q}_{k-1}$ for any $\gamma$, we get $$\UT_{k-1}(\Delta_k(\gamma)) \, \widehat{q}_{k-1} \, = \, \sum_{i \le j \le \ell} \pi_{i,k-1} d\gamma_k \pi_{j,k-1} \widehat{q}_{k-1}.$$ Therefore, if we know that $d\gamma_k = A_k + B_k$ where the left support of $A_k$ and the right support of $B_k$ are dominated by $\sum_{s \ge 1} p_{k-s} = \1_\A - q_{k-1}$, then we deduce that $\UT_{k-1}(\Delta_k(\gamma)) \, \widehat{q}_{k-1} = 0$. In other words, it suffices to prove that $$q_{k-1} \Delta_k(\gamma) q_{k-1} \, = \, 0  \quad \mbox{for} \quad \gamma = b_d, g_{\mathit{off}}, b_{\mathit{off}}.$$ We have 
\begin{eqnarray*}
\Delta_k(b_d) & = & \sum_j \Delta_k \big( p_j(f-f_j)p_j \big) \\ & = & \sum_{j < k} p_j (f_k - f_j) p_j - \sum_{j < k-1} p_j (f_{k-1} - f_j) p_j \\ & = & \sum_{j \le k-1} p_j df_k p_j \ = \ (\1_\A - q_{k-1}) \Delta_k(b_d) (\1_\A - q_{k-1}).
\end{eqnarray*}
To calculate the martingale differences for $g_{\mathit{off}}$, we invoke the formula $$g_{\mathit{off}} = \sum_{s=1}^\infty \sum_{j \in \Z} p_j df_{j+s} q_{j+s-1} + q_{j+s-1} df_{j+s} p_j$$ given in the statement of the Calder\'on-Zygmund decomposition. Then we find 
\begin{eqnarray*}
\Delta_k(g_{\mathit{off}}) & = & \sum_{s=1}^\infty p_{k-s} df_k q_{k-1} + q_{k-1} df_k p_{k-s} \\ & = & (\1_\A - q_{k-1}) df_k q_{k-1} + q_{k-1} df_k (\1_\A - q_{k-1}).
\end{eqnarray*}
Finally, it remains to consider the martingale differences of $b_{\mathit{off}}$
\begin{eqnarray*}
\Delta_k(b_{\mathit{off}}) & = & \sum_{s=1}^\infty \sum_{j \in \Z} \Delta_k \big( p_j (f - f_{j+s}) p_{j+s} + p_{j+s} (f - f_{j+s}) p_j \big) \\ & = & \sum_{s=1}^\infty \, \sum_{j < k-s} p_j (f_k - f_{j+s}) p_{j+s} + p_{j+s} (f_k - f_{j+s}) p_j \\ & - & \sum_{s=1}^\infty \sum_{j < k-s-1} p_j (f_{k-1} - f_{j+s}) p_{j+s} + p_{j+s} (f_{k-1} - f_{j+s}) p_j \\ & = & \sum_{s=1}^\infty \, \sum_{j < k-s} p_j df_k p_{j+s} + \sum_{s=1}^\infty \, \sum_{j < k-s} p_{j+s} df_k p_j \ = \ A_k + B_k.
\end{eqnarray*}
So $q_{k-1} A_k \hskip-1pt = \hskip-1pt B_k q_{k-1} \hskip-1pt = \hskip-1pt 0$ and $q_{k-1} \Delta_k(\gamma) q_{k-1} \hskip-1pt = \hskip-1pt 0$ for $\gamma \hskip-1pt = \hskip-1pt b_d, g_{\mathit{off}}, b_{\mathit{off}}$ as desired. \fin

\subsection{Haar shift operators} 

The Haar system has the form $$h_Q^\varepsilon(x) \, = \, \frac{1}{\sqrt{|Q|}} \prod_{j=1}^n \big( 1_{I_j^-}(x_j) + \varepsilon_j 1_{I_j^+}(x_j) \big)$$ where $Q = I_1 \times I_2 \times \cdots \times I_n \in \Q$ and $\varepsilon = (\varepsilon_1, \varepsilon_2, \ldots, \varepsilon_n) \neq (1,1,\ldots,1)$ with $\varepsilon_j \in \pm 1$. We are using $I_j^-$ and $I_j^+$ for the left/right halves of the intervals $I_j$. It yields an orthonormal system in $L_2(\R^n)$ composed of mean zero functions. If we write $h_Q$ for any Haar function of the form $h_Q^\varepsilon$, a \emph{noncommuting dyadic shift with complexity $(r,s)$} has the form $$\Sha_\alpha f(x) \, = \, \sum_{Q \in \Q} A_Qf \, = \, \sum_{Q \in \Q} \ \sum_{\begin{subarray}{c} R,S \ \mathrm{dyadic} \, \subset Q \\ \ell(R) = 2^{-r} \ell(Q) \\ \ell(S) = 2^{-s} \ell(Q)\end{subarray}} \alpha_{RS}^Q \big\langle f, h_S \big\rangle h_R(x),$$ where $\langle f, h_S \rangle = \int f h_S$ and $\alpha_{RS}^Q$ are operators in $\M$ satisfying $\|\alpha_{RS}^Q\|_\M \le \frac{\sqrt{|R||S|}}{|Q|}$.

\begin{lemma} \label{L2}
We have $\|\Sha_\alpha f \|_2 \le \|f\|_2$.
\end{lemma}

\dem The argument is standard, observe that 
\begin{eqnarray*}
\|\Sha_\alpha f\|_2^2 & = & \sum_{Q,Q'} \sum_{R,R'S,S'} \tau \big( \langle f, h_S \rangle^* \alpha_{RS}^{Q*} \alpha_{R'S'}^{Q'} \langle f, h_{S'} \rangle \big) \, \int_{\R^n} h_R(y) h_{R'}(y) \, dy.
\end{eqnarray*}
The integral on the right imposes $R=R'$, which in turn gives $Q=Q'$ since $Q$ is the unique $r$-th ancestor of $R$ and the same happens for $(R',Q')$. Once we know that $Q=Q'$, we may write $$\|\Sha_\alpha f\|_2^2 \, = \, \sum_{Q \in \Q} \|A_Q f\|_2^2 \, = \, \sum_{Q \in \Q} \Big\| A_Q \Big( \sum_{\begin{subarray}{c} S \subset Q \\ \ell(S) = 2^{-s} \ell(Q) \end{subarray}} \langle f, h_S \rangle h_S \Big) \Big\|_2^2.$$ It is worth mentioning that the double use above of $h_S$ always refers to the same choice of $h_S^\varepsilon$ in both instances. On the other hand, it is easily seen that $A_Q$ is a contractive operator on $L_2(\A)$. Indeed, we have 
\begin{eqnarray*}
\|A_Qg\|_2^2 & \le & \int_{\R^n} \Big[ \sum_{R,S} \|\alpha_{R,S}^Q\|_\M \Big( \frac{1}{\sqrt{|S|}} \int_S \|g(y)\|_{L_2(\M)} dy \Big) \frac{1}{\sqrt{|R|}} 1_R(x) \Big]^2 \, dx \\ & \le & \int_Q \Big( \mean_Q \|g(y)\|_{L_2(\M)} \, dy \Big)^2 \, dx \ \le \ \mean_Q \|g\|_{L_2(\A)}^2 \, dx \ = \ \|g\|_{L_2(\A)}^2.
\end{eqnarray*}
This yields \\ [8pt] \null \hskip5pt \hfill $\displaystyle \|\Sha_\alpha f\|_2^2 \, \le \, \sum_{Q \in \Q} \Big\| \sum_{\begin{subarray}{c} S \subset Q \\ \ell(S) = 2^{-s} \ell(Q) \end{subarray}} \langle f, h_S \rangle h_S \Big\|_2^2 \ = \ \Big\| \sum_{Q \in \Q} \langle f, h_Q \rangle h_Q \Big\|_2^2 \ = \ \|f\|_2^2.$ \hfill \fin

The next lemma is crucial to analyze Haar shifts and general Calder\'on-Zygmund operators with noncommuting kernels. We take here the opportunity to slightly modify the argument in \cite[Lemma 4.2]{Pa1}, which was not entirely correct.  

\begin{lemma} \label{Dilation}
Given $s \in \Z_+$, there exists $\zeta \in \A_\pi$ such that$\, :$
\begin{itemize}
\item[i)] $\lambda \varphi (\1_\A - \zeta) \, \le \, 2^{sn} \|f\|_1$,

\vskip3pt

\item[ii)] If $Q_0 \in \Q_{k_0}$ and $x \in \widehat{Q}_0^s$, then $\zeta(x) \le \widehat{q}_{k_0}(y)$ for all $y \in Q_0$. 
\end{itemize}
In the second property, we write $\widehat{Q}_0^s$ for the unique $s$-th dyadic ancestor of $Q_0$.
\end{lemma}

\dem We have $$\1_\A - \widehat{q}_k \, = \, \sum_{j \le k} \big( \widehat{q}_{j-1} - \widehat{q}_j \big) \, = \, \sum_{j \le k} \sum_{Q \in \Q_j} \rho_Q \otimes 1_Q \, = \, \sum_{Q \in \Q_k} \Big[ \sum_{R \supset Q} \rho_R \Big] \otimes 1_Q$$ for some family of projections $\rho_Q \in \M_\pi$. Define $$\zeta \, = \, \bigwedge_{k \in \Z} \zeta_k \quad \mbox{with} \quad \zeta_k \, = \, \1_\A - \bigvee_{j \le k} \bigvee_{Q \in \Q_j} \rho_Q 1_{\widehat{Q}^s}.$$ It is clear that the $\zeta_k$'s are decreasing in $k$ and we find
\begin{eqnarray*}
\lambda \varphi(\1_\A - \zeta) & = & \lambda \lim_{k \to \infty} \varphi (\1_\A - \zeta_k) \\ & \le & \lambda \lim_{k \to \infty} \sum_{j \le k} \sum_{Q \in \Q_j} \tau(\rho_Q) |\widehat{Q}^s| \\ & = & 2^{sn} \lim_{k \to \infty} \lambda \sum_{j \le k} \sum_{Q \in \Q_j} \varphi(\rho_Q \otimes 1_Q) \\ & = & 2^{sn} \lambda \, \varphi \big( \1_\A - \widehat{q} \hskip1pt \big) \ = \ 2^{sn} \lambda \sum_{m \ge \ell} \varphi \big( \1_\A - q(2^m) \big)  \ \lesssim \ 2^{sn} \|f\|_1.
\end{eqnarray*}
To prove the second property, it will be useful to observe that $Q_1 \subsetneq Q_2$ implies that $\rho_{Q_1} \perp \rho_{Q_2}$ are orthogonal projections. Indeed, according to the definition of $\rho_Q$ above, we have $\rho_{Q_1} \rho_{Q_2} 1_{Q_1} = (\widehat{q}_{j_1-1} - \widehat{q}_{j_1}) (\widehat{q}_{j_2-1} - \widehat{q}_{j_2}) 1_{Q_1} = 0$ for $\ell(Q_1) = 2^{-j_1}$ and $\ell(Q_2) = 2^{-j_2}$. Then, we find
\begin{eqnarray*}
\hskip50pt \zeta(x) & \le & \zeta_{k_0}(x) \\ [10pt] & = & \1_\M - \bigvee_{j \le k_0} \bigvee_{Q \in \Q_j} \rho_Q 1_{\widehat{Q}^s}(x) \\ & \le & \1_\M - \bigvee_{R \supset Q_0} \rho_R \ = \ 1_\M - \sum_{R \supset Q_0} \rho_R \\ & = & \Big( \1_\A - \sum_{Q \in \Q_{k_0}} \Big[ \sum_{R \supset Q} \rho_R \Big] \otimes 1_Q \Big) (y) \ = \ \widehat{q}_{k_0}(y). \hskip60pt \square
\end{eqnarray*}

\demAib As in the perfect dyadic case, we assume $f \in \A_{c,+}$ and decompose $f = f_r + f_c$ in the same way. Once more the argument is row/column symmetric, and we just consider the column part. After fixing $\lambda = 2^\ell$ for some $\ell \in \Z$, we construct the corresponding Calder\'on-Zygmund decomposition for $f_c = g_d^c + g_{\mathit{off}}^c + b_d^c + b_{\mathit{off}}^c$. According to Lemma \ref{L2}, we may control the term $\Sha_\alpha(g_d^c)$ in the usual way. Given $\gamma \in \{b_d, g_{\mathit{off}}, b_{\mathit{off}} \}$, the other terms can be decomposed as follows
\begin{eqnarray*}
\Sha_\alpha (\gamma^c) & = & \sum_{k \in \Z} \Sha_\alpha \big( \mathsf{UT}_{k-1}(\Delta_k(\gamma)) \big) \\ [15pt] & = & \sum_{k \in \Z} \hskip1pt \sum_{Q \in \Q} \hskip1pt \sum_{\begin{subarray}{c} R,S \subset Q \\ \ell(R) = 2^{-r} \ell(Q) \\ \ell(S) = 2^{-s} \ell(Q) \end{subarray}} \alpha_{RS}^Q \Big( \int_{\R^n} \mathsf{UT}_{k-1}(\Delta_k(\gamma)) h_S \, dy \Big) h_R(x) \\ & = & \sum_{k \in \Z} \hskip1pt \Big[ \sum_{\begin{subarray}{c} Q \in \Q \\ \ell(Q) \le 2^{-k+1} \end{subarray}} \hskip1pt \ + \ \sum_{\begin{subarray}{c} Q \in \Q \\ \ell(Q) > 2^{-k+1} \\ \ell(Q) \le 2^{s-k+1} \end{subarray}} \ + \ \sum_{\begin{subarray}{c} Q \in \Q \\ \ell(Q) > 2^{s-k+1} \end{subarray}} \Big] \ = \ A_\gamma + B_\gamma + C_\gamma.
\end{eqnarray*}
We claim that $C_\gamma = 0$. Namely, we have $\ell(S) = 2^{-s} \ell(Q) > 2^{-k+1}$. This means that $\mathsf{E}_{k-1}(h_S) = h_S$ since the Haar functions $h_S$ are constant in the dyadic children of $S$, whose length sides are greater or equal than $2^{-(k-1)}$. This yields 
\begin{eqnarray*}
\int_{\R^n} \mathsf{UT}_{k-1}(\Delta_k(\gamma)) h_S \, dy  & = & \int_{\R^n} \mathsf{E}_{k-1} \big( \mathsf{UT}_{k-1}(\Delta_k(\gamma)) h_S \big) \, dy \\ & = & \int_{\R^n} \mathsf{E}_{k-1} \big( \mathsf{UT}_{k-1}(\Delta_k(\gamma)) \big) h_S \, dy \\ & = & \int_{\R^n} \big( \mathsf{UT}_{k-1}( \mathsf{E}_{k-1} \Delta_k(\gamma)) h_S \, dy \ = \ 0.
\end{eqnarray*}
In order to deal with the remaining terms $A_\gamma$ and $B_\gamma$, we invoke the identity $q_{k-1} \Delta_k(\gamma) q_{k-1} = 0$ which was already justified in the perfect dyadic case whenever $\gamma = b_d, g_{\mathit{off}}, b_{\mathit{off}}$. Namely, since $\pi_{i,k-1} (\1_\A - q_{k-1}) = (\1_\A - q_{k-1}) \pi_{j,k-1} = 0$ for $i,j \le \ell$, we find $$\mathsf{UT}_{k-1}(\Delta_k(\gamma)) \, = \, \sum_{i \le j} \pi_{i,k-1} \Delta_k(\gamma) \pi_{j,k-1} \, = \, \sum_{\begin{subarray}{c} i \le j \\ j > \ell \end{subarray}} \pi_{i,k-1} \Delta_k(\gamma) \pi_{j,k-1}.$$ Let us now consider the term $A_\gamma$, we have $$\lambda \, \varphi \big\{ |A_\gamma| > \lambda \big\} \, \le \, \lambda \, \varphi \big( \1_\A - \widehat{q} \hskip1pt \big) \, + \,  \lambda \, \varphi \Big\{ \big| A_\gamma \widehat{q} \hskip1pt \big| > \frac{\lambda}{2} \Big\}.$$ We already know that the first term on the right is dominated by $\|f\|_1$ and 
\begin{eqnarray*}
A_\gamma \widehat{q} & = & \sum_{k \in \Z} \hskip1pt \sum_{\begin{subarray}{c} Q \in \Q \\ \ell(Q) \le 2^{-k+1} \end{subarray}} \hskip1pt \sum_{\begin{subarray}{c} R,S \subset Q \\ \ell(R) = 2^{-r} \ell(Q) \\ \ell(S) = 2^{-s} \ell(Q) \end{subarray}} \alpha_{RS}^Q \Big( \int_{\R^n} \mathsf{UT}_{k-1}(\Delta_k(\gamma)) h_S \, dy \Big) h_R(x) \, \widehat{q}(x).
\end{eqnarray*} 
Given $Q \in \Q$ with $\ell(Q) \le 2^{-k+1}$ let $$k_Q \ge k-1 \quad \mbox{determined by} \quad \ell(Q) = 2^{-k_Q}.$$ It is clear that $\widehat{q}(x) = \widehat{q}_{k_Q}(x) \widehat{q}(x) = \widehat{q}_{k_Q}(y) \widehat{q}(x) = \widehat{q}_{k-1}(y) \widehat{q}(x)$ whenever $x,y$ belong to $Q$. However, the presence of $h_R(x), h_S(y)$ implies (unless the corresponding term is $0$) that the pair $(x,y) \in R \times S \subset Q \times Q$ so that we may write $$A_\gamma \widehat{q} \, = \, \sum_{k \in \Z} \hskip1pt \sum_{\begin{subarray}{c} Q \in \Q \\ \ell(Q) \le 2^{-k+1} \end{subarray}} \hskip1pt \sum_{\begin{subarray}{c} R,S \subset Q \\ \ell(R) = 2^{-r} \ell(Q) \\ \ell(S) = 2^{-s} \ell(Q) \end{subarray}} \alpha_{RS}^Q \Big( \int_{\R^n} \mathsf{UT}_{k-1}(\Delta_k(\gamma)) \widehat{q}_{k-1} h_S \, dy \Big) h_R(x) \, \widehat{q}(x).$$ Therefore, we conclude $$\mathsf{UT}_{k-1}(\Delta_k(\gamma)) \widehat{q}_{k-1} \, = \, \sum_{\begin{subarray}{c} i \le j \\ j > \ell \end{subarray}} \pi_{i,k-1} \Delta_k(\gamma) \pi_{j,k-1} \widehat{q}_{k-1} \, = \, 0$$ since $\pi_{j,k-1} \widehat{q}_{k-1} = 0$ when $j > \ell$. This shows that $A_\gamma \widehat{q} = 0$. Let us finally consider the term $B_\gamma$. We will follow a similar argument with the projection $\zeta$ from Lemma \ref{Dilation} instead. Namely, we have $$\lambda \, \varphi \big\{ |B_\gamma| > \lambda \big\} \, \le \, \lambda \, \varphi \big( \1_\A - \zeta \big) \, + \,  \lambda \, \varphi \Big\{ \big| B_\gamma \zeta \hskip1pt \big| > \frac{\lambda}{2} \Big\}.$$ According to property i) of Lemma \ref{Dilation}, it suffices to show that $B_\gamma \zeta = 0$. Now we know that $\ell(Q) \le 2^{s-k+1}$, so that $k_Q \ge k-s-1$. Let us now consider the $2^{ns}$ dyadic cubes $T_j$ having $Q$ as their $s$-th dyadic ancestor. This gives rise to the identities $$\zeta(x) \, = \, \zeta_{k_Q + s}(x) \zeta(x) \, = \, \zeta_{k_Q + s}(y) \zeta(x) \, = \, \widehat{q}_{k_Q + s}(z) \zeta(x) \ = \ \widehat{q}_{k-1}(z) \zeta(x)$$ for $(x,y,z) \in Q \times Q \times T_j$. Indeed, the second identity follows from the fact that $\mathsf{E}_{k_Q}(\zeta_{k_Q+s}) = \zeta_{k_Q+s}$, the third one from the second property in Lemma \ref{Dilation} and the last one from the inequality $k_Q \ge k-s-1$. Hence, given $y \in S \subset Q$ we pick the unique $j$ for which $S = T_j$ and deduce that $\zeta(x) = \widehat{q}_{k-1}(y) \zeta(x)$. Then it yields the identity $$B_\gamma \zeta \, = \, \sum_{k \in \Z} \hskip1pt \sum_{\begin{subarray}{c} Q \in \Q \\ \ell(Q) > 2^{-k+1} \\ \ell(Q) \le 2^{s-k+1} \end{subarray}} \hskip1pt \sum_{\begin{subarray}{c} R,S \subset Q \\ \ell(R) = 2^{-r} \ell(Q) \\ \ell(S) = 2^{-s} \ell(Q) \end{subarray}} \alpha_{RS}^Q \Big( \int_{\R^n} \mathsf{UT}_{k-1}(\Delta_k(\gamma)) \widehat{q}_{k-1} h_S \, dy \Big) h_R(x) \, \zeta(x).$$ The integrand $\mathsf{UT}_{k-1}(\Delta_k(\gamma)) \widehat{q}_{k-1}$ vanishes for the same reason as it did above. 
\fin

\begin{remark}
\emph{Our constants are $\sim 2^{sn}$ and seem far to be sharp. Unfortunately, the classical argument leading to constants $\sim s$ encounters a major obstacle due to the presence ---in the noncommutative setting--- of triangular truncations, which are not bounded in $L_1$. The Appendix below contains more details on this topic.}
\end{remark}

\subsection{Noncommuting CZO's} 

The proofs of Theorems Aii), B and C arise from a careful combination of recent results in the theory of noncommutative Hardy spaces. Let us begin introducing Mei's notion \cite{MeiMAMS} of row and column Hardy spaces for our algebra of operator-valued functions $\A$. In order to distinguish from order Hardy spaces to be introduced below, let us follows Mei's notation and define $$\mathrm{H}_1(\R^n; \M) \, = \, \mathrm{H}_1^r(\R^n; \M) + \mathrm{H}_1^c(\R^n;\M)$$ as the space of functions $f \in L_1(\A)$ for which we have $$\|f\|_{\mathrm{H}_1(\R^n;\M)} \, = \, \inf_{f = g+h} \|g\|_{\mathrm{H}_1^r(\R^n;\M)} + \|h\|_{\mathrm{H}_1^c(\R^n;\M)} < \infty,$$ where the row/column norms are given  by 
\begin{eqnarray*}
\|g\|_{\mathrm{H}_1^r(\R^n;\M)} & = & \Big\| \Big( \int_\Gamma \Big[ \frac{\partial \widehat{g}}{\partial t} \hskip0.5pt \frac{\partial \widehat{g}^*}{\partial t} \hskip0.5pt + \hskip0.5pt \summ_j \frac{\partial \widehat{g}}{\partial x_j} \hskip0.5pt \frac{\partial \widehat{g}^*}{\partial x_j} \Big] (x + \cdot,t) \, \frac{dx dt}{t^{n-1}} \Big)^\frac12 \Big\|_1, \\ [-1pt] \|h\|_{\mathrm{H}_1^c(\R^n;\M)} & = & \Big\| \Big( \int_\Gamma \Big[ \frac{\partial \widehat{h}^*}{\partial t} \frac{\partial \widehat{h}}{\partial t} + \summ_j \frac{\partial \widehat{h}^*}{\partial x_j} \frac{\partial \widehat{h}}{\partial x_j} \Big] (x + \cdot,t) \, \frac{dx dt}{t^{n-1}} \Big)^\frac12 \Big\|_1,
\end{eqnarray*}
with $\Gamma = \{ (x,t) \in \R^{n+1}_+ \, | \ |x| < y \}$ and $\widehat{f}(x,t) = P_t f(x)$ for the Poisson semigroup $(P_t)_{t \ge 0}$. In other words, operator-valued forms of Lusin's square function. We say that $a \in L_1(\M; L_2^c(\R^n))$ is a \emph{column atom} if there exists a cube $Q$ so that 
\begin{itemize}
\item $\mathrm{supp}_{\R^n} \hskip1pt a = Q$, 

\vskip8pt

\item $\displaystyle \int_Q a(y) \, dy = 0$,

\item $\|a\|_{L_1(\M;L_2^c(\R^n))} = \displaystyle \tau \Big[ \big( \int_Q |a(y)|^2 \, dy \big)^\frac12 \Big] \le \frac{1}{\sqrt{|Q|}}$.
\end{itemize}
According to \cite[Theorem 2.8]{MeiMAMS}, we have $$\|f\|_{\mathrm{H}_1^c(\R^n;\M)} \, \sim \, \inf \Big\{ \summ_k |\lambda_k| \, \big| \ f = \summ_k \lambda_k a_k \ \mbox{with} \ a_k \ \mbox{column atoms} \Big\}.$$ On the other hand, we have already settled a dyadic filtration $(\A_k)_{k \in \Z}$ for our algebra of operator-valued functions $\A$. Then, we may follow \cite{PX1} to define the corresponding noncommutative Hardy space $\mathrm{H}_1(\A)$ as the completion of the space of finite martingales in $L_1(\A)$ with respect to the norm $$\|f\|_{\mathrm{H}_1(\A)} \, = \, \inf_{\begin{subarray}{c} f = g + h \\ g,h \ \mathrm{martingales} \end{subarray}} \Big\| \Big( \sum_{k \in \Z} dg_k dg_k^* \Big)^\frac12 \Big\|_1 + \Big\| \Big( \sum_{k \in \Z} dh_k^* dh_k \Big)^\frac12 \Big\|_1.$$ In other words, $\mathrm{H}_1(\A) = \mathrm{H}_1^r(\A) + \mathrm{H}_1^c(\A)$ where the spaces on the right are the completions of the spaces of finite $L_1$-martingales with respect to the norms in $L_1$ of the corresponding row/column square functions given above.  By the use of a dyadic covering \cite{CGP,MeiMAMS}, it can be shown that there exists $n+1$ dyadic filtrations $\Sigma_\A^j$ $(0 \le j \le n)$ in $\R^n$ so that $$\mathrm{H}_1(\R^n;\M) \, \simeq \, \sum_{j=0}^n \mathrm{H}_1(\A, \Sigma_\A^j),$$ where the latter spaces are defined as $\mathrm{H}_1(\A)$ after replacing the standard filtration $\Sigma_\A^0$ by any other dyadic filtration in our family. Moreover, this isomorphism also holds independently for row/column Hardy spaces. 

\vskip3pt

\demAii It suffices to show $$T_r: \mathrm{H}_1^r(\A) \to L_1(\A) \quad \mbox{and} \quad T_c: \mathrm{H}_1^c(\A) \to L_1(\A),$$ for any generic noncommuting CZO $(T_r,T_c)$. Indeed, in that case we decompose $f = f_r + f_c \in \mathrm{H}_1(\A)$, so that $\|f\|_{\mathrm{H}_1(\A)} \sim \|f_r\|_{\mathrm{H}_1^r(\A)} + \|f_c\|_{\mathrm{H}_1^c(\A)}$ and we deduce that $$\|T_rf_r\|_1 + \|T_cf_c\|_1 \, \lesssim \, \|f_r\|_{\mathrm{H}_1^r(\A)} + \|f_c\|_{\mathrm{H}_1^c(\A)} \sim \, \|f\|_{\mathrm{H}_1(\A)}.$$ According to our observation above, $\mathrm{H}_1(\A)$ embeds isomorphically into $\mathrm{H}_1(\R^n;\M)$ by means of a suitably choice of dyadic coverings of $\R^n$, and the same holds for row and column spaces isolatedly. Therefore, it also suffices to show that $$T_r: \mathrm{H}_1^r(\R^n;\M) \to L_1(\A),$$ $$T_c: \mathrm{H}_1^c(\R^n;\M) \to L_1(\A).$$ Both estimates are identical, let us prove the column case. According to the atomic decomposition of $\mathrm{H}_1^c(\R^n;\M)$ we just find a uniform upper estimate for the $L_1$ norm of $T_c(a)$ valid for an arbitrary column atom $$\|T_c(a)\|_1 \, \le \, \big\| T_c(a) 1_{2Q} \big\|_1 + \big\| T_c(a) 1_{\R^n \setminus 2Q} \big\|_1.$$ The second term is dominated by  
\begin{eqnarray*}
\big\| T_c(a) 1_{\R^n \setminus 2Q} \big\|_1 & = & \tau \int_{\R^n \setminus 2Q} \Big| \int_Q k(x,y) a(y) \, dy \Big| \, dx \\ & \le & \int_Q \Big( \int_{\R^n \setminus 2Q} \big\| k(x,y) - k(x,c_Q) \big\|_\M \, dx \Big) \tau |a(y)| \, dy \\ [6pt] & \lesssim & \tau \Big( \int_Q |a(y)| \, dy \Big) \ \le \ \sqrt{|Q|} \tau \Big[ \big( \int_Q |a(y)|^2 \, dy \big)^\frac12 \Big] \ \le \ 1,
\end{eqnarray*}
where the next to last estimate follows from Hansen's inequality or as a consequence of the operator-convexity of the function $a \mapsto | a |^2$. As for the first term, it suffices to show that $T_c: L_1(\M; L_2^c(\R^n)) \to L_1(\M; L_2^c(\R^n))$, since then we find again 
\begin{eqnarray*}
\big\| T_c(a) 1_{2Q} \big\|_1 & = & \tau \Big( \int_{2Q} |T_c(a)(x)| \, dx \Big) \\ & \le & \sqrt{|2Q|} \, \tau \Big[ \big( \int_{2Q} |T_c(a)(x)|^2 \, dx \big)^\frac12 \Big] \\ & \lesssim & \sqrt{|2Q|} \, \tau \Big[ \big( \int_{Q} |a(x)|^2 \, dx \big)^\frac12 \Big] \, \lesssim \, 1.
\end{eqnarray*}
The $L_1(\M; L_2^c(\R^n))$-boundedness of $T_c$ follows from anti-linear duality $$\big\| T_c(f) \big\|_{L_1(\M;L_2^c(\R^n))} \, \le \, \Big( \sup_{\|g\|_{L_\infty(L_2^c)} \le 1} \big\|T_c^*(g) \big\|_{L_\infty(\M;L_2^c(\R^n))} \Big) \|f\|_{L_1(\M;L_2^c(\R^n))}.$$ It is easily checked that the adjoint $T_c^*(g)$ has the form $T_c^*g(x) \sim \int_{\R^n} k(y,x)^* g(y) \, dy$ when we construct it with respect to the anti-linear bracket $\langle f,g \rangle = \varphi (f* g)$. This means in particular that $T_c^*$ is still an $L_2$-bounded column CZO associated to a kernel satisfying H\"ormander smoothness. This gives rise to 
\begin{eqnarray*}
\big\| T_c^*(g) \big\|_{L_\infty(\M;L_2^c(\R^n))} & = & \Big\| \Big( \int_{\R^n} |T_c^*(g)(x)|^2 \, dx \Big)^\frac12 \Big\|_\M \\ & = & \sup_{\|u\|_{L_2(\M)} \le 1} \Big( \int_{\R^n} \big\langle |T_c^*(g)(x)|^2 u,u \big\rangle_{L_2(\M)} \, dx \Big)^\frac12 \\ & = & \sup_{\|u\|_{L_2(\M)} \le 1} \Big( \int_{\R^n} \big\|T_c^*(gu)(x) \big\|_{L_2(\M)}^2 \, dx \Big)^\frac12 \\ & \lesssim & \sup_{\|u\|_{L_2(\M)} \le 1} \Big( \int_{\R^n} \big\| g(x)u \big\|_{L_2(\M)}^2 \, dx \Big)^\frac12 \\ & = & \Big\| \Big( \int_{\R^n} |g(x)|^2 \, dx \Big)^\frac12 \Big\|_\M.
\end{eqnarray*}
The third identity above uses the right $\M$-module nature of column CZO's. \fin

\begin{remark}
\emph{Theorem Aii) could have also been derived from the $L_\infty \to \mathrm{BMO}$ type estimates in \cite{JMP1}. We have preferred to include this alternative argument using atomic decompositions. Still a third approach is possible using more recent atomic decompositions from \cite{Pe2,HM}. This will be needed below for martingale transforms and paraproducts. The proof goes in fact a little further than the statement, since it emphasizes row/column $\mathrm{H}_1 \to L_1$ type estimates for $T_r/T_c$ respectively. This also works for arbitrary semicommutative CZO's under suitable assumptions, see \cite{JMP1} for details.}
\end{remark}

\begin{remark} \label{Rem-Modularity}
\emph{The proof above also shows that $L_1(L_2^\dag)$ and $L_\infty(L_2^\dag)$ boundedness of $T_\dag$ for $\dag \in \{r,c\}$ follow from the corresponding $L_2$ boundedness of the same operator. As noticed in \cite{JMP1}, this is very specific of CZO's with noncommuting kernels since other semicommutative CZO's fail to satisfy this implication. The key property here is left/right $\M$-modularity, so that $$uT_r(f) = T_r(uf) \quad \mbox{and} \quad T_c(f)u = T_c(fu).$$ This also explains our approach through weak type estimates, see the Appendix.}
\end{remark}

\subsection{Row/column $L_p$ estimates}

Theorem B follows as an easy consequence of Theorem A after applying suitable interpolation/duality results. Thus, we will only outline the definition of the involved spaces and the necessary results to deduce Theorem B from Theorem A. Given $1 < p < \infty$, the noncommutative Hardy space $\mathrm{H}_p(\A)$ is defined as $$\mathrm{H}_p(\A) \, = \, \begin{cases} \mathrm{H}_p^r(\A) + \mathrm{H}_p^c(\A) & \mbox{if} \ 1 < p \le 2, \\ \mathrm{H}_p^r(\A) \cap \mathrm{H}_p^c(\A) & \mbox{if} \ 2 \le p < \infty, \end{cases}$$ where the corresponding row/column Hardy spaces arise as the completion of the subspace of finite martingales in $L_p(\A)$ with respect to the norms given by the row and column square functions 
\begin{eqnarray*}
\|f\|_{\mathrm{H}_p^r(\A)} & = & \Big\| \Big( \sum_{k \in \Z} df_k df_k^* \Big)^\frac12 \Big\|_p, \\ \|f\|_{\mathrm{H}_p^c(\A)} & = & \Big\| \Big( \sum_{k \in \Z} df_k^* df_k \Big)^\frac12 \Big\|_p.
\end{eqnarray*}
Pisier/Xu obtained in \cite{PX1} the noncommutative Burkholder-Gundy inequalities which can be formulated as $L_p(\A) \simeq \mathrm{H}_p(\A)$ for $1 < p < \infty$. On the other hand, we know from \cite{J1,JX} that $\mathrm{H}_p^\dag(\A)^* \simeq \mathrm{H}_{p'}^\dag(\A)$ for $\dag \in \{r,c\}$ and $1 < p < \infty$. Regarding interpolation, we know from Musat \cite{Musat} that $$\mathrm{H}_p^\dag(\A) \, \simeq \, \big[ \mathrm{H}_{p_0}^\dag(\A), \mathrm{H}_{p_1}^\dag(\A) \big]_\theta,$$ where $\dag \in \{r,c\}$ and $\frac1p = \frac{1-\theta}{p_0} + \frac{\theta}{p_1}$. The proof of Theorem B is now straightforward.

\vskip3pt

\demB We know that $$T_r: \mathrm{H}_1^r(\A) \to L_1(\A) \quad \mbox{and} \quad T_c: \mathrm{H}_1^c(\A) \to L_1(\A).$$ If $1 < p < 2$, we find $T_r: \mathrm{H}_p^r(\A) \to L_p(\A)$ and $T_c: \mathrm{H}_p^c(\A) \to L_p(\A)$ by interpolation with $L_2(\A) = \mathrm{H}_2^r(\A) = \mathrm{H}_2^c(\A)$. Hence, taking a decomposition $f = f_r + f_c$ satisfying $\|f\|_p \sim \|f\|_{\mathrm{H}_p(\A)} \sim \|f_r\|_{\mathrm{H}_p^r(\A)} + \|f_c\|_{\mathrm{H}_p^c(\A)}$ we get $\|T_rf_r \|_p + \| T_cf_c \|_p \lesssim \|f\|_p$. Now if $2 < p < \infty$, recalling that $T_r^*, T_c^*$ are again row/column CZO's with the same properties, duality gives $T_r: L_p(\A) \to \mathrm{H}_p^r(\A)$ and $T_c: L_p(\A) \to \mathrm{H}_p^c(\A)$. This immediately yields the inequality in Theorem Bii). The $L_\infty \to \mathrm{BMO}$ type estimates were originally proved in \cite{JMP1}, these also follows by duality from Theorem A. \fin

\begin{remark}
\emph{Alternatively, it can be proved that the row/column $L_p$ estimates in Theorem Bi) for $1 < p < 2$ also follow by real interpolation from the weak type estimates in Theorem Ai). Moreover, since Mei's spaces $\mathrm{H}_p(\R^n;\M)$ also behave well for interpolation and duality, the statement of Theorem B could have been done in terms of these other Hardy spaces.}
\end{remark}

\section{Proof of Theorem C}

In this section we turn our attention to noncommutative martingale transforms and paraproducts. In particular, the former pair $(\A,\varphi)$ will refer in what follows to an arbitrary semifinite von Neumann algebra equipped with a normal faithful semifinite trace. Our filtration $\Sigma_\A = (\A_k)_{k \ge 1}$ will be any increasing family of von Neumann subalgebras, whose union is weak-$*$ dense in $\A$. The operators $\mathsf{E}_k$ and $\Delta_k$ still denote the corresponding conditional expectations and martingale difference operators. As mentioned in the Introduction, we will deal with
\begin{itemize} 
\item[\textbf{a)}] \textbf{Noncommuting martingale transforms} $$M_\xi^r f = \sum_{k \ge 1} \Delta_k(f) \xi_{k-1} \quad \mbox{and} \quad M_\xi^c f = \sum_{k \ge 1} \xi_{k-1} \Delta_k(f).$$ 

\vskip3pt

\item[\textbf{b)}] \textbf{Paraproducts with noncommuting symbol} $$\Pi_\rho^r(f) = \sum_{k \ge 1} \mathsf{E}_{k-1}(f) \Delta_k(\rho)  \quad \mbox{and} \quad \Pi_\rho^c(f) = \sum_{k \ge 1} \Delta_k(\rho) \mathsf{E}_{k-1}(f).$$
\end{itemize}
The martingale coefficients $\xi_{k} \in \A_{k}$ form an adapted sequence and it is easy to show that $L_2$-boundedness of $M_\xi^r$ and $M_\xi^c$ holds iff the $\xi_k$'s are uniformly bounded in the norm of $\A$. On the other hand, the classical characterization $\Pi_\rho: L_2 \to L_2$ iff $\rho \in \mathrm{BMO}$ was disproved by Nazarov, Pisier, Treil and Volberg \cite{NPTV}, see also Mei's paper \cite{M1}. Hence, the $L_2$-boundedness of $\Pi_\rho^r$ and $\Pi_\rho^c$ will be simply assumed in what follows. Regarding Cuculescu's construction and CZ decomposition, no essential changes are needed. Namely, given $f \in L_1^+(\A)$ (the former space $\A_{c,+}$ is unnecessary since our filtration starts now at $k=1$) and $\lambda \in \R_+$, Cuculescu's construction is verbatim the same. The only difference is on the diagonal estimate $$\Big\| qfq + \sum_{k=1}^\infty p_k f_k p_k \Big\|_2^2 \, \lesssim \, \lambda \|f\|_1.$$ This inequality requires to work with regular filtrations, which are defined through the additional condition $\mathsf{E}_k(f) \le c \hskip1pt \mathsf{E}_{k-1}(f)$ for some absolute constant $c>0$ and every pair $(f,k) \in \A_+ \times \Z_+$. Of course, the reader might think that it is more appropriate to use in this case the noncommutative form of Gundy's decomposition \cite{PR}, which does not require any regularity assumption on the martingale. This leads unfortunately to some problems related to our triangular truncations which will be explained in the Appendix below.   

\vskip3pt

\demCi The argument is essentially the same as in the perfect dyadic case. Given $f \in L_1^+(\A)$, we construct the same decomposition $f = f_r + f_c$ via the projections $\pi_{j,k}$ and fix $\lambda = 2^\ell$ for some $\ell \in \Z$. A further CZ decomposition gives $f_c = g_d^c + g_{\mathit{off}}^c + b_d^c + b_{\mathit{off}}^c$ as usual. According to our regularity assumption, we still have $$\max \Big\{ \|g_d^r\|_2^2, \|g_d^c\|_2^2 \Big\} \, \le \, \|g_d\|_2^2 \, = \, \Big\| qfq + \sum_{k \ge 1} p_kf_kp_k \Big\|_2^2 \, \lesssim \, \lambda \|f\|_1.$$ Thus, arguing as in the proof of Theorem A it suffices to show that $$\widehat{q} \hskip1pt M_\xi^r(\gamma^r) \, = \, M_\xi^c(\gamma^c) \hskip1pt \widehat{q} \, = \, \widehat{q} \hskip1pt \Pi_\rho^r(\gamma^r) \, = \, \Pi_\rho^c(\gamma^c) \hskip1pt \widehat{q} \, = \, 0$$ for any $\gamma \in \{g_{\mathit{off}}, b_d, b_{\mathit{off}} \}$. As usual, we just consider the column case by symmetry. Let us begin with martingale transforms. Since $\gamma^c = \sum_j \mathsf{UT}_{j-1}(\Delta_j(\gamma))$ and the triangular truncation $\mathsf{UT}_{j-1}$ is built with $j$-predictable projections, we see that $\mathsf{UT}_{j-1}(\Delta_j(\gamma))$ is a $j$-th martingale difference, so that $$\Delta_k (\gamma^c) \, = \, \mathsf{UT}_{k-1}(\Delta_k(\gamma)).$$ By the proof of Theorem A, we know $\mathsf{UT}_{k-1}(\Delta_k(\gamma)) \hskip1pt \widehat{q}_{k-1} = 0$ and $$M_\xi^c(\gamma^c) \hskip1pt \widehat{q} \, = \, \sum_{k = 1}^\infty \xi_{k-1} \Delta_k(\gamma^c) \hskip1pt \widehat{q} \, = \, \sum_{k = 1}^\infty \xi_{k-1} \mathsf{UT}_{k-1}(\Delta_k(\gamma)) \hskip1pt \widehat{q}_{k-1} \hskip1pt \widehat{q} \, = \, 0.$$ For martingale paraproducts, we observe that $\mathsf{E}_{k-1}(\gamma^c) = \sum_{j < k} \mathsf{UT}_{j-1}(\Delta_j(\gamma))$ and \\ \null \hskip2.5cm  $\displaystyle \Pi_\rho^c(\gamma^c) \hskip1pt \widehat{q} \, = \, \sum_{k=1}^\infty \Delta_k(\rho) \sum_{j < k} \mathsf{UT}_{j-1}(\Delta_j(\gamma)) \hskip1pt \widehat{q}_{j-1} \hskip1pt \widehat{q} \, = \, 0.$ \hfill \fin

\begin{remark}
\emph{Adjoints of martingale paraproducts have the form $$\big[ \Pi_\rho^c \big]^* f \, = \, \sum_{k \ge 1} \mathsf{E}_{k-1} \big( \Delta_k(\rho^*) \Delta_k(f) \big) \quad \mbox{and} \quad \big[ \Pi_\rho^r \big]^* f \, = \, \sum_{k \ge 1} \mathsf{E}_{k-1} \big( \Delta_k(f) \Delta_k(\rho^*) \big)$$ when using the anti-linear duality bracket. It is easy to adapt the argument above for these maps, to obtain weak type inequalities for adjoints of noncommutative paraproducts associated to regular filtrations $$\inf_{f = f_r + f_c} \big\| \big[ \Pi_\rho^r \big]^*f_r \big\|_{1,\infty} + \big\| \big[ \Pi_\rho^c \big]^*f_c \big\|_{1,\infty} \, \le \, \|f\|_1.$$} 
\end{remark}

We defined above the noncommutative Hardy spaces $\mathrm{H}_1(\A)$. Alternatively, we may also consider the noncommutative form $\mathrm{h}_1(\A) \, = \, \mathrm{h}_1^r(\A) + \mathrm{h}_1^c(\A) + \mathrm{h}_1^d(\A)$ of the conditional Hardy space $\mathrm{h}_1$, where the norms are given by 
\begin{eqnarray*}
\|f\|_{\mathrm{h}_1^r(\A)} & = & \Big\| \Big( \sum_{k \ge 1} \mathsf{E}_{k-1} \big( df_k df_k^* \big) \Big)^\frac12 \Big\|_1, \\ \|f\|_{\mathrm{h}_1^c(\A)} & = & \Big\| \Big( \sum_{k \ge 1} \mathsf{E}_{k-1} \big( df_k^* df_k \big) \Big)^\frac12 \Big\|_1,\\ [5pt] \|f\|_{\mathrm{h}_1^d(\A)} & = & \Big\| \sum_{k \ge 1} | df_k | \Big\|_1 \ = \ \sum_{k \ge 1} \|df_k\|_1.
\end{eqnarray*}
The space $\mathrm{h}_1(\A)$ was studied in \cite{JM,Pe1}, it was independently proved that 
\begin{eqnarray*}
\mathrm{H}_1^r(\A) & \simeq & \mathrm{h}_1^r(\A) + \mathrm{h}_1^d(\A), \\
\mathrm{H}_1^c(\A) & \simeq & \mathrm{h}_1^c(\A) + \mathrm{h}_1^d(\A).
\end{eqnarray*}
In conjunction, these isomorphisms could be regarded as a noncommutative form of Davis' decomposition for martingales. Shortly after, it was found in \cite{Pe2} an atomic decomposition for the spaces $\mathrm{h}_1^r(\A)$ and $\mathrm{h}_1^c(\A)$. More precisely, an element $a$ in $L_1(\A) \cap L_2(\A)$ is called a \emph{column atom} with respect to the filtration $(\A_k)_{k \ge 1}$ if there exists $k_0 \in \Z_+$ and a finite projection $e \in \A_{k_0}$ such that
\begin{itemize}
\item $a = ae$,

\vskip3pt

\item $\mathsf{E}_{k_0}(a) = 0$,

\vskip2pt

\item $\|a\|_2 \le \varphi(e)^{-\frac12}$.
\end{itemize}
An element $a \in L_1(\A)$ is called a $\mathrm{c-atom}$ if it is a column atom or $a \in \A_1$ with $\|a\|_1 \le 1$. Row atoms are defined to satisfy $a=ea$ instead and $\mathrm{r-atoms}$ are defined similarly. We also refer to \cite{HM} for $q$-analogs of these notions. In the following result, we collect some norm equivalences coming from atomic decompositions and John-Nirenberg type inequalities. Recall that 
\begin{eqnarray*}
\|f\|_{\mathrm{BMO}_c(\A)} & = & \sup_{k \ge 1} \Big\| \mathsf{E}_k \big[ (f - f_{k-1})^*(f - f_{k-1}) \big] \Big\|_\A^\frac12, \\ \|f\|_{\mathrm{bmo}_c(\A)} \hskip4pt & = & \max \Big\{ \big\| \mathsf{E}_1(f) \big\|_1, \ \sup_{k \ge 1} \Big\| \mathsf{E}_k \big[ (f-f_k)^*(f-f_k) \big] \Big\|_\A^\frac12 \Big\}. 
\end{eqnarray*}
As usual, the corresponding row norms of $f$ arise as the column norms of $f^*$. If we also define $\|f\|_{\mathrm{bmo}_d(\A)} = \sup_k \|df_k\|_\A$, then we can define the spaces $\mathrm{BMO}(\A)$ and $\mathrm{bmo}(\A)$ as follows
\begin{eqnarray*}
\|f\|_{\mathrm{BMO}(\A)} & = & \max \Big\{ \|f\|_{\mathrm{BMO}_r(\A)}, \|f\|_{\mathrm{BMO}_c(\A)} \Big\}, \\ \|f\|_{\mathrm{bmo}(\A)} \hskip4pt & = & \max \Big\{ \|f\|_{\mathrm{bmo}_r(\A)}, \|f\|_{\mathrm{bmo}_c(\A)}, \|f\|_{\mathrm{bmo}_d(\A)} \Big\}.  
\end{eqnarray*}
The isomorphism $\mathrm{BMO}(\A) \simeq \mathrm{bmo}(\A)$ was independently proved in \cite{JM,Pe1}. 

\begin{Atomic} We have 
\begin{eqnarray*}
\|f\|_{\mathrm{h}_1^r} & \sim & \inf \Big\{ \summ_k |\lambda_k| \, \big| \ f = \summ_k \lambda_k a_k \ \mathrm{and} \ a_k \ \mathrm{r-atom} \Big\}, \\ [3pt] \|f\|_{\mathrm{h}_1^c} & \sim & \inf \Big\{ \summ_k |\lambda_k| \, \big| \ f = \summ_k \lambda_k a_k \ \mathrm{and} \ a_k \ \mathrm{c-atom} \Big\}, \\ [3pt] \|f\|_{\mathrm{bmo}(\A)} & \sim &  \sup_{k \ge 1} \Big[ \|df_k\|_\infty \vee \sup_{\begin{subarray}{c} \beta \in \A_k \\ \|\beta\|_1 \le 1 \end{subarray}} \big\| \beta (f-f_k) \big\|_1 \vee \sup_{\begin{subarray}{c} \beta \in \A_k \\ \|\beta\|_1 \le 1 \end{subarray}} \big\| (f-f_k) \beta \big\|_1\Big].
\end{eqnarray*}
The last equivalence is a John-Nirenberg type inequality, which differs from \cite{JMus}.
\end{Atomic} 

\demCii Let us begin with $\mathrm{H}_1 \to L_1$ type inequalities. As pointed out in the proof of Theorem Aii), it suffices to show that $T_\dag: \mathrm{H}_1^\dag(\A) \to L_1(\A)$ with $\dag \in \{r,c\}$ and for both martingale transforms and paraproducts. Since we have $$\mathrm{H}_1^\dag(\A) \simeq \mathrm{h}_1^\dag(\A) + \mathrm{h}_1^d(\A),$$ it suffices to show that $T_\dag: \mathrm{X} \to L_1(\A)$ with $\mathrm{X}$ any of the two spaces appearing on the right. Once more, the argument is row/column symmetric and we just consider columns. To see that $T_c: \mathrm{h}_1^c(\A) \to L_1(\A)$ we may use the atomic decomposition above, so that it suffices to find a uniform upper bound for $\|T_c(a)\|_1$ with $a$ being a $\mathrm{c-atom}$. If $a \in \A_1$ with $\|a\|_1 \le 1$, then we see that $$M_\xi^c(a) = \xi_0 a_1 \quad \mbox{and} \quad \Pi_\rho^c(a) = ba = \Pi_\rho^c(u|a|^\frac12) |a|^\frac12 \quad \mbox{for} \quad a = u|a|.$$ In particular, $\|M_\xi^c(a)\|_1 + \|\Pi_\rho^c(a)\|_1 \lesssim \|a\|_1 \le 1$. If $a$ is a column atom, we find 
\begin{eqnarray*}
M_\xi^c(a) & = & \sum_{k > k_0} \xi_{k-1} \Delta_k(a) \ = \ \sum_{k > k_0} \xi_{k-1} \Delta_k(a) e \ = \ M_\xi^c(a)e, \\ \Pi_\rho^c(a) & = & \sum_{k > k_0+1} \Delta_k(\rho) \mathsf{E}_{k-1}(a) \ = \ \sum_{k > k_0+1} \Delta_k(\rho) \mathsf{E}_{k-1}(a) e \ = \ \Pi_\rho^c(a)e.
\end{eqnarray*}
This gives rise to $\|T_c(a)\|_1 = \|T_c(a)e\|_1 \le \|T_c(a)\|_2 \|e\|_2 \lesssim \|a\|_2 \|e\|_2 \le 1$ for both martingale transforms and paraproducts. We have already justified the $\mathrm{h}_1^c \to L_1$ boundedness. Let us now look at $\mathrm{h}_1^d$ $$\|M_\xi^c(f)\|_1 \, \le \, \sum_{k \ge 1} \|\xi_k\|_\infty \|\Delta_k(f)\|_1 \, \le \, \Big( \sup_{k \ge 1} \|\xi_k\|_\infty \Big) \hskip1pt \|f\|_{\mathrm{h}_1^d(\A)}$$ As for the paraproduct, we use the John-Nirenberg inequality above 
\begin{eqnarray*}
\|\Pi_\rho^c(f)\|_1 & = & \Big\| \sum_{k \ge 1} \Delta_k(\rho) \sum_{j < k} \Delta_j(f) \Big\|_1 \\ & = & \Big\| \sum_{k \ge 1} \big( \rho - \rho_k \big) \Delta_k(f) \Big\|_1 \ \lesssim \ \|\rho\|_{\mathrm{bmo}(\A)} \|f\|_{\mathrm{h}_1^d(\A)}.
\end{eqnarray*}
According to \cite{JM,Pe1} and \cite{M1,NPTV}, we have $$\|\rho\|_{\mathrm{bmo}(\A)} \, \sim \, \|\rho\|_{\mathrm{BMO}(\A)} \, \lesssim \, \max \Big\{ \big\| \Pi_\rho^r: L_2 \to L_2 \big\|, \big\| \Pi_\rho^c: L_2 \to L_2 \big\| \Big\}.$$ All together gives that $M_\xi^c$ and $\Pi_\rho^c$ take $\mathrm{H}_1^c(\A)$ into $L_1(\A)$ as we claimed. In fact slight modifications of the given argument yield the same result for $[\Pi_\rho^c]^*$, details are left to he reader. This is all what is needed to produce analog inequalities in this setting to those in Theorems A and B, we just need to follow the arguments verbatim. It remains to show that $\Pi_\rho^c: L_p(\A) \to L_p(\A)$ for $p > 2$, for which it will be enough to prove $L_\infty \to \mathrm{BMO}$ boundedness and use interpolation. The $L_\infty \to \mathrm{BMO}_c$ boundedness follows by duality from the $\mathrm{H}_1^c \to L_1$ boundedness of $[\Pi_\rho^c]^*$. On the other hand, the $L_\infty \to \mathrm{BMO}_r$ boundedness is very simple 
\begin{eqnarray*}
\|\Pi_\rho^cf\|_{\mathrm{BMO}_r(\A)} & = & \sup_{k \ge 1} \Big\| \mathsf{E}_k \Big( \sum_{j \ge k} \Delta_j(\Pi_\rho^c(f)) \Delta_j(\Pi_\rho^c(f))^* \Big) \Big\|_\A^\frac12 \\ & = & \sup_{k \ge 1} \Big\| \mathsf{E}_k \Big( \sum_{j \ge k} \Delta_j(\rho) \mathsf{E}_{j-1}(f) \mathsf{E}_{j-1}(f)^* \Delta_j(\rho)^* \Big) \Big\|_\A^\frac12 \\ & \le & \sup_{k \ge 1} \Big\| \mathsf{E}_k \Big( \sum_{j \ge k} \Delta_j(\rho) \Delta_j(\rho)^* \Big) \Big\|_\A^\frac12 \, \|f\|_\infty \ \le \ \|\rho\|_{\mathrm{BMO}_r(\A)} \|f\|_\infty.
\end{eqnarray*} 
Now we majorize $\|\rho\|_{\mathrm{BMO}_r(\A)}$ by the $L_2 \to L_2$ norm of $\Pi_\rho$ as we did above. \fin

Observe that we have not needed to assume regularity of our martingale filtration and we find that $[\Pi_\rho^r]^*, [\Pi_\rho^c]^*$ take $\mathrm{H}_1 \to L_1$ and $L_p \to L_p$ for $1 < p < 2$ by duality. In some sense, row/column noncommutative paraproducts present a similar behavior as row/column square functions in the noncommutative  Burkholder-Gundy and Khintchine inequalities \cite{Lu,LuP,PX1}. On the other hand, \cite[Theorem 5.7]{Rad2} yields $L \log L \to L_1$ type estimates for a finite von Neumann algebra $\A$ with $(T_r,T_c)$ a martingale transform/paraproduct with noncommuting coefficients/symbol $$\inf_{f = f_r + f_c} \big\| T_rf_r \big\|_1 + \big\| T_cf_c \big\|_1 \, \lesssim \, \|f\|_{L \log L(\A)}.$$

\section*{Appendix. Open problems} 

\subsection*{A.1. CZO's with noncommuting kernels}

Our proof of Theorem Ai) is not entirely satisfactory, since it does not include arbitrary CZO's with noncommuting kernels. In the general case, we can not expect to annihilate the terms associated to $g_{\mathit{off}}, b_d, b_{\mathit{off}}$. If the reader considers the simplest term $b_d$, a difficulty with triangular truncations in $L_1$ will be immediately recognized. In fact, our proof for Haar shifts operators does not provide sharp constants for the same reason.

\vskip3pt

\noindent \textbf{Problem 1.} Extend Theorem Ai) to arbitrary CZO's with noncommuting kernels. 

\vskip3pt

Here is a possible alternative argument. Once we have $f = f_r + f_c$, the same decomposition constructed in the proof of the perfect dyadic case, we could consider a \emph{left \emph{CZ} decomposition} for $f_r$ and a \emph{right \emph{CZ} decomposition} for $f_c$ as follows. Given $\lambda \in \R_+$ we let $f_r = g_r + b_r$ and $f_c = g_c + b_c$ with $$g_r \, = \, \widehat{q} \hskip1pt f_r + \sum_{k \in \Z} \widehat{p}_k \mathsf{E}_k(f_r) \quad \mbox{and} \quad b_r \, = \, \sum_{k \in \Z} \widehat{p}_k \big( f_r - \mathsf{E}_k(f_r) \big),$$ where $\widehat{p}_k = \widehat{q}_{k-1} - \widehat{q}_k$. The column decomposition just requires to put $\widehat{p}_k$ and $\widehat{q}$ on the right. The advantage of this approach is that we do not find off-diagonal terms which were much harder to deal in \cite{Pa1}. Moreover, it is not very difficult to show that $$\max \Big\{ \|g_r\|_2^2, \|g_c\|_2^2 \Big\} \, \lesssim \, \lambda \|f\|_1$$ as expected. Problem 1 would be solved if we knew that $$\sum_{k \in \Z} \big\| \widehat{p}_k \big( f_r - \mathsf{E}_k(f_r) \big) \big\|_1 + \big\| \big( f_c - \mathsf{E}_k(f_c) \big) \widehat{p}_k \big\|_1 \, \lesssim \, \|f\|_1.$$ It is perhaps too optimistic to expect that the inequality above holds, since the triangular truncations $\mathsf{LT}_k$ and $\mathsf{UT}_k$ appear to be incomparable for different values of $k$. We wonder whether some noncommutative form of Davis' decomposition in the sense of \cite{Rad3} could be useful to modify our row/column decomposition $f = f_r + f_c$ before performing the CZ decomposition, see also \cite{Pa0} for related ideas. Note that such a row/column CZ decomposition would provide in particular a much simpler proof of the main result in \cite{Pa1}, since off diagonal terms would disappear. 

\vskip3pt

\noindent \textbf{Problem 2.} Find a row/column CZ decomposition of $f$ in the line explained above.  

\subsection*{A.2. CZO's on general von Neumann algebras}

As explained in \cite{Pa1}, a key ingredient for a successful application of the noncommutative CZ decomposition is to use it on $\M$-bimoludar maps. In this paper, our decomposition $f = f_r + f_c$ has allowed us to make it work for either left or right $\M$-module maps. There are however many other semicommutative CZO's, some of which were mentioned in the Introduction. We know from \cite{JMP1} that a semicommutative CZO satisfying $L_\infty(L_2^r)$ and $L_\infty(L_2^c)$ boundedness also satisfies $T: L_\infty(\A) \to \mathrm{BMO}(\A)$.

\vskip3pt

\noindent \textbf{Problem 3.} Do we have $T: L_1(\A) \to L_{1,\infty}(\A)$ under the same assumptions?

\vskip3pt

According to \cite{JMP1}, solving Problem 3 for CZO's associated to a kernel acting by Schur multiplication would provide weak type $(1,1)$ inequalities for crossed product extensions of classical CZO's $$Tf(x) \, \sim \, \sum_{g \in \G} \int_{\R^n} k(x,y) f_g(y) \rtimes_\gamma \lambda(g) \, dy$$ on $\A = L_\infty(\R^n) \rtimes_\gamma \G$. This in turn is closely related to weak type estimates for Fourier multipliers on group von Neumann algebras, see \cite{JMP1} for further details. On the other hand, consider CZO's of the form $$Tf(x) \, \sim \, \int_{\R^n} (id \otimes \mathrm{tr}) \Big[ k(x,y) \big( \mathbf{1} \otimes f(y) \big) \Big] \, dy.$$ As we have seen along this paper and in \cite{Pa1}, weak type inequalities require to find vanishing products $q_1(y) q_2(x)$ with $q_1, q_2$ certain projections in $\A$, see e.g. Lemma \ref{Dilation}. However, we find $T(fq_1)(x) q_2(x) \sim \int_{\R^n} (id \otimes \mathrm{tr}) [ k(x,y) ( q_2(x) \otimes fq_1(y) ) ] dy$ in the model above and no interaction between $q_1$ and $q_2$ takes place. This is due to the lack of right $\M$-modularity for $T$. In fact, solving Problem 3 for this kind of CZO's is very much related to the CZ theory for von Neumann algebras developed in \cite{JMP3}. Namely, the projection in Lemma \ref{Dilation} is a dyadic dilation on $\R^n$ of $\widehat{q}$ not affecting its $\M$ \lq structure\rq${}$ because the CZO is given as a partial trace on $\R^n$, but not on $\M$. The idea in the model above is to dilate both in $\R^n$ and $\M$. Dilating in $\M$ has to do with finding a suitable \lq metric\rq${}$ in $\M$ to work with. This is what is done in \cite{JMP3} in terms of diffusion semigroups on the given algebra. Under this point of view, we could relate CZO's on $(\A, \varphi)$ with those in \cite{Pa1} when $\varphi$ is tracial and with the ones considered in this paper when $\varphi$ is a nontracial weight. 

\vskip3pt

\noindent \textbf{Problem 4.} Prove a CZ decomposition/weak type inequalities for CZO's in \cite{JMP3}.

\subsection*{A.3. Gundy's decomposition vs triangular truncations}

It is a little bit unsatisfactory to require regular filtrations to provide weak type inequalities for martingales transforms/paraproducts with noncommuting coefficients/symbols. It is well-known that these estimates hold in the classical setting for any filtration by means of Gundy's decomposition. The noncommutative extension of Gundy's decomposition was constructed in \cite{PR}. Given a positive martingale $f = (f_1, f_2, \ldots)$ in $L_1(\A)$, we may decompose it as $f = \alpha + \beta + \gamma$ with 
\begin{eqnarray*}
d\alpha_k & = & q_k df_k q_k - \mathsf{E}_{k-1} \big( q_k df_k q_k \big), \\
d\beta_k & = & q_{k-1}df_kq_{k-1} - q_k df_kq_k + \mathsf{E}_{k-1} \big( q_k df_k q_k \big), \\
d\gamma_k & = & df_k - q_{k-1} df_k q_{k-1}.
\end{eqnarray*}
It was proved in \cite{PR} that $$\max \Big\{ \frac{1}{\lambda} \|\alpha\|_2^2, \sum_{k \ge 1} \|d\beta_k\|_1, \lambda \varphi \Big( \bigvee_{k \ge 1} \mathrm{supp}^* d\gamma_k \Big) \Big\} \, \lesssim \, \|f\|_1,$$ where $\mathrm{supp}^* a = \mathbf{1}_\A - q$ with $q$ the greatest projection satisfying $q a q = 0$. If we try to prove Theorem C using Gundy's decomposition instead of Calder\'on-Zygmund decomposition, we will not find any trouble controlling the terms associated to $\alpha$ and $\gamma$. The term $\beta$ presents however a significant difficulty due to the presence of triangular truncations $\mathsf{LT}_k$ and $\mathsf{UT}_k$ in $L_1(\A)$. This difficulty can be summarized as follows. Consider a family $\mathsf{Tr}_k$ of upper triangular truncations and assume that $(\alpha_k, \beta_k) \in L_\infty(\A) \times L_1(\A)$, do we have $$\Big\| \sum_{k \ge 1} \alpha_k \mathsf{Tr}_k(\beta_k) \Big\|_1 \, \lesssim \, \sum_{k \ge 1} \|\alpha_k \beta_k\|_1\mbox{?}$$ Or at least $$\Big\| \sum_{k \ge 1} \alpha_k \mathsf{Tr}_k(\beta_k) \Big\|_1 \, \lesssim \, \Big( \sup_{k \ge 1} \|\alpha_k\|_\infty \Big) \, \sum_{k \ge 1} \|\beta_k\|_1\mbox{?}$$ The first condition suffices to manage paraproducts with noncommuting symbols, the second one is weaker but sufficient to deal with martingale transforms having noncommuting coefficients. When dealing with lower triangular truncations, we should have $\mathsf{Tr}_k(\beta_k) \alpha_k$ on the left and $\beta_k \alpha_k$ on the right hand side. 

\vskip3pt

\noindent \textbf{Problem 5.} Does any of these inequalities hold?

\vskip3pt

\noindent \textbf{Problem 6.} Can we eliminate the regularity assumption from Theorem C?

\vskip3pt

\noindent \textbf{Acknowledgement.} We would like to thank Tao Mei for some discussions we had with him on the content of this paper. Guixiang Hong was supported in part by an ANR Grant 2011 BS01 008 01 (France). Luis Daniel L\'opez-S\'anchez and Jos\'e Mar\'ia Martell were supported in part by MEC Grant MTM-2010-16518 (Spain). Javier Parcet was supported in part by the ERC Grant StG-256997 (European Union).

\bibliographystyle{amsplain}

\vskip30pt

\hfill \noindent \textbf{Guixiang Hong} \\
\null \hfill Laboratoire de Math{\'e}matiques
\\ \null \hfill Universit{\'e} de France-Comt{\'e} \\
\null \hfill 16 Route de Gray, 25030 Besan\c{c}on Cedex, France \\
\null \hfill\texttt{ghong@univ-fcomte.fr}

\

\hfill \noindent \textbf{Luis Daniel L\'opez} \\
\null \hfill Instituto de Ciencias Matem{\'a}ticas \\ \null \hfill
CSIC-UAM-UC3M-UCM \\ \null \hfill Consejo Superior de
Investigaciones Cient{\'\i}ficas \\ \null \hfill C/ Nicol\'as Cabrera 13-15.
28049, Madrid. Spain \\ \null \hfill\texttt{luisd.lopez@icmat.es}

\

\hfill \noindent \textbf{Jos\'e Mar\'ia Martell} \\
\null \hfill Instituto de Ciencias Matem{\'a}ticas \\ \null \hfill
CSIC-UAM-UC3M-UCM \\ \null \hfill Consejo Superior de
Investigaciones Cient{\'\i}ficas \\ \null \hfill C/ Nicol\'as Cabrera 13-15.
28049, Madrid. Spain \\ \null \hfill\texttt{chema.martell@icmat.es}

\

\hfill \noindent \textbf{Javier Parcet} \\
\null \hfill Instituto de Ciencias Matem{\'a}ticas \\ \null \hfill
CSIC-UAM-UC3M-UCM \\ \null \hfill Consejo Superior de
Investigaciones Cient{\'\i}ficas \\ \null \hfill C/ Nicol\'as Cabrera 13-15.
28049, Madrid. Spain \\ \null \hfill\texttt{javier.parcet@icmat.es}
\end{document}